\documentclass[11pt]{article}
\usepackage{amsmath}
\usepackage{dsfont}
\usepackage{mathrsfs}
\usepackage{amsmath,amssymb}
\usepackage{amsfonts}
\usepackage{hyperref}
\usepackage{amsthm}
\usepackage{graphicx}
\usepackage{subfigure}
\usepackage{xcolor}
\renewcommand{\qed}{\hfill\small{$\blacksquare$}\normalsize}

\hfuzz=\maxdimen
\theoremstyle{definition}
\newtheorem{lemma}{Lemma}[section]
\newtheorem{definition}[lemma]{Definition}
\newtheorem{proposition}[lemma]{Proposition}
\newtheorem{theorem}[lemma]{Theorem}
\newtheorem{corollary}[lemma]{Corollary}
\newtheorem{example}{Example}
\newtheorem{remark}{Remark}
\newtheorem{question}{Question}

\newtheorem{conjecture}{Conjecture}
\numberwithin{equation}{section}
\renewcommand{\proof}{\textbf{Proof. }}
\renewcommand{\qed}{\hfill\small{$\blacksquare$}\normalsize}

\DeclareFixedFont{\Acknowledgment}{OT1}{cmr}{bx}{n}{14pt}
\begin{document}

\title{\bf A combinatorial Yamabe problem on two and three dimensional manifolds}
\author{Huabin Ge, Xu Xu}
\maketitle

\begin{abstract}
In this paper, we define a new discrete curvature on two and three dimensional triangulated manifolds,
which is a modification of the well-known discrete curvature on these manifolds.
The new definition is more natural and respects the scaling exactly the same way as Gauss curvature does.
Moreover, the new discrete curvature can be used to approximate the Gauss curvature on surfaces.
Then we study the corresponding constant curvature problem, which is called the combinatorial Yamabe problem,
by the corresponding combinatorial versions of Ricci flow and Calabi flow for surfaces
and Yamabe flow for 3-dimensional manifolds.
The basic tools are the discrete maximal principle and variational principle.
\end{abstract}

\tableofcontents

\section{Introduction}\label{Introduction}
\subsection{Background and Preliminaries}\label{Introduction-Background}
One of the central topics in differential geometry is the existence of canonical metrics on a given smooth manifold,
especially metrics with constant curvatures.
One can ask similar questions on triangulated manifolds.
Different from the standard processes on smooth manifolds, where all calculations are done in local coordinates,
Regge \cite{Re} suggested an entirely different way, where the calculations are done in a geometric simplex.
Since the geometric information in a single simplex is included in the lengths of edges,
the geometric approach to a triangulated manifold may be started with a piecewise flat metric,
which assigns each edge a Euclidean distance such that each simplex can be realized as a Euclidean simplex.
Piecewise flat metric brings singularities on codimension two sub-simplices,
which are explained as simplicial curvature tensors in \cite{Al,CMS,MMAGWY,MMAGY}.
Besides defining edge lengths directly, one can also define discrete metrics at vertices and derive the edge lengths indirectly. In fact, Thurston \cite{T1} introduced the notion of circle packing metric on a triangulated surface, while Cooper and Rivin introduced the notion of sphere packing metric on a triangulated 3-manifold. We shall study the circle (sphere) packing metrics and the corresponding discrete curvatures in this paper, and we want to know when there exist canonical discrete metrics and how to find them.

Suppose $M$ is a closed surface with a triangulation $\mathcal{T}=\{V,E,F\}$,
where $V,E,F$ represent the sets of vertices, edges and faces respectively.
Let $\Phi: E\rightarrow [0,\frac{\pi}{2}]$ be a function assigning each edge $\{ij\}$ a weight $\Phi_{ij}\in [0,\frac{\pi}{2}]$.
The triple $(M, \mathcal{T}, \Phi)$ will be referred to as a weighted triangulation of $M$ in the following.
All the vertices are ordered one by one, marked by $v_1, \cdots, v_N$, where $N=V^\sharp$
is the number of vertices. We use $i\sim j$ to denote that the
vertices $i$ and $j$ are adjacent if there is an edge $\{ij\}\in E$ with $i$, $j$ as end points.
Throughout this paper, all functions $f: V\rightarrow \mathds{R}$ will be regarded as column
vectors in $\mathds{R}^N$ and $f_i$ is the value of $f$ at $i$. And we use $C(V)$ to denote the set of functions
defined on $V$. Each map $r:V\rightarrow (0,+\infty)$ is called a circle packing metric.
Given $(M, \mathcal{T}, \Phi)$, we attach each edge $\{ij\}$ a length
\begin{equation}\label{definition of length of edge}
l_{ij}=\sqrt{r_i^2+r_j^2+2r_ir_j\cos \Phi_{ij}}.
\end{equation}
Thurston proved \cite{T1} that the lengths $\{l_{ij}, l_{jk}, l_{ik}\}$ satisfy the triangle inequalities,
which ensures that the face $\{ijk\}$ could be
realized as an Euclidean triangle $\triangle v_iv_jv_k$ with lengths $\{l_{ij}, l_{jk}, l_{ik}\}$.
In this sense, the triangulated surface $(M, \mathcal{T}, \Phi)$ could be taken as
gluing many Euclidean triangles coherently. However, this gluing procedure produces singularities at the vertices,
which are described by the discrete curvatures.
Suppose $\theta_i^{jk}$ is the inner angle of the triangle $\triangle v_iv_jv_k$
at the vertex $i$, the classical well-known discrete Gaussian curvature at $i$ is defined as
\begin{equation}\label{classical Gauss curv}
K_i=2\pi-\sum_{\triangle v_iv_jv_k \in F}\theta_i^{jk},
\end{equation}
where the sum is taken over all the triangles with $v_i$ as one of its vertices. Discrete Gaussian curvature $K_i$ satisfies the following discrete version of Gauss-Bonnet formula \cite{CL1}:
\begin{equation}\label{Gauss-Bonnet without weight}
\sum_{i\in V}K_i=2\pi \chi(M).
\end{equation}
Denote $K_{av}=\frac{2\pi \chi(M)}{N}$, it is natural to ask the following combinatorial Yamabe problem:

\begin{question}\label{Question-2d-K-curvature}
Given a weighted triangulated surface $(M, \mathcal{T}, \Phi)$, does there exist a circle packing metric $r$
that determines a constant $K$-curvature? How to find it if it exists?
\end{question}

Different from the smooth case, the constant $K$-curvature metric, i.e. a metric $r$ with $K_i=K_{av}$ for all $i\in V$, does not always exist. Thurstion \cite{T1} found that, besides the topological structure, the \textbf{combinatorial structure} plays an essential role for the existence of constant $K$-curvature metric. In fact, for any proper subset $I\subset V$, let $F_I$ be the subcomplex whose vertices are in $I$ and let $Lk(I)$ be the set of pairs $(e, v)$ of an edge $e$ and a vertex $v$ satisfying the following three conditioins: (1) The end points of $e$ are not in $I$; (2) $v$ is in $I$; (3) $e$ and $v$ form a triangle. Thurston \cite{T1} proved that

\begin{theorem}\label{Thm-thurston} \textbf{(Thurston)} \;Given a weighted triangulated surface $(M, \mathcal{T}, \Phi)$,
the existence of constant $K$-curvature metric is equivalent to the following combinatorial-topological conditions
\begin{equation}\label{condition-Thurston}
2\pi\chi(M)\frac{|I|}{|V|} >-\sum_{(e,v)\in Lk(I)}(\pi-\Phi(e))+2\pi\chi(F_I), \;\;\forall I: \phi\subsetneqq I\subsetneqq V.
\end{equation}
Moreover, the constant $K$-curvature metric is unique, if exists, up to the scaling of $r$.
\end{theorem}

In fact, the following result is proved for the classical discrete curvature $K$ \cite{T1,MR,CL1}.

\begin{theorem}\label{Thm-thurston-MardenRodin-ChowLuo} \textbf{(Thurston-Marden\&Rodin-Chow\&Luo)}
\;Given a weighted triangulated surface $(M, \mathcal{T}, \Phi)$, denote
\begin{equation}\label{Def-Y}
\mathscr{Y}\triangleq\mathscr{K}_{GB} \cap (\mathop{\cap} _{\phi\neq I \subsetneqq V} \mathscr{Y}_I ),
\end{equation}
where
\begin{equation}\label{Gauss-Bonnet plane}
\mathscr{K}_{GB}\triangleq\Big\{x\in \mathds{R}^N\Big|\sum_{i=1}^Nx_i=2\pi\chi(M)\Big\}
\end{equation}
and for any nonempty proper subset $I\subset V$,
\begin{equation}\label{Y_I half space}
\mathscr{Y}_I\triangleq\Big\{x\in \mathds{R}^N \Big|\sum_{i\in I}x_i >-\sum_{(e,v)\in Lk(I)}(\pi-\Phi(e))+2\pi\chi(F_I)\Big\}.
\end{equation}
The space of all admissible $K$-curvatures $\big\{K=K(r)\big|r\in \mathds{R}^N_{>0}\big\}$ is exactly $\mathscr{Y}$.
\end{theorem}

Chow and Luo \cite{CL1} first established an intrinsic connection between Thurston's circle packing metric and the surface Ricci flow.
They introduced a combinatorial Ricci flow
\begin{equation*}\label{Def-ChowLuo's flow}
\frac{dr_i}{dt}=-K_ir_i
\end{equation*}
with the normalization
\begin{equation}\label{Def-ChowLuo's normalized flow}
\frac{dr_i}{dt}=(K_{av}-K_i)r_i.
\end{equation}
Then they reproved Theorem \ref{Thm-thurston} by the curvature flow methods. In fact, they proved:

\begin{theorem}\label{Thm-Chow-Luo} \textbf{(Chow\&Luo)} \;Given a weighted triangulated surface $(M, \mathcal{T}, \Phi)$,
the solution to the normalized combinatorial Ricci flow (\ref{Def-ChowLuo's normalized flow}) converges (exponentially fast to a constant $K$-curvature metric) if and only if there exists a constant $K$-curvature metric.
\end{theorem}

Inspired by \cite{CL1}, the first author introduced a combinatorial Calabi flow
\begin{equation}\label{Def-Ge-Calabi-flow}
\frac{dr_i}{dt}=\Delta K_ir_i
\end{equation}
in \cite{Ge} and proved similar results:
\begin{theorem}\label{Thm-Ge-Calabi} \textbf{(Ge)} \;
Given a weighted triangulated surface $(M, \mathcal{T}, \Phi)$,
the solution to the combinatorial Calabi flow (\ref{Def-Ge-Calabi-flow}) converges (exponentially fast to a constant $K$-curvature metric)
if and only if there exists a constant $K$-curvature metric.
\end{theorem}

Theorem \ref{Thm-thurston}, Theorem \ref{Thm-thurston-MardenRodin-ChowLuo}, Theorem \ref{Thm-Chow-Luo} and Theorem \ref{Thm-Ge-Calabi} answer Question \ref{Question-2d-K-curvature} perfectly.
The first two theorems give a necessary and sufficient combinatorial-topological condition for the existence of a constant $K$-curvature metric,
while the last two theorems give two different efficient ways to find a constant $K$-curvature metric, i.e.
for any initial value $r(0)$, the solutions to the combinatorial Ricci flow and Calabi flow converge automatically to a constant $K$-curvature metric exponentially fast.

\subsection{Main results}\label{Introduction-main-results}
It seems that the classical discrete Gauss curvature $K$ is a suitable analogue for smooth Gauss curvature.
However, there are some intrinsic disadvantages for the classical discrete Gauss curvature.
For one thing, the smooth Gauss curvature changes under the scaling of Riemannian metrics
while the discrete Gauss curvature dose not.
For another, the discrete Gauss curvature $K$ can not be used directly to approximate the smooth Gauss curvature
when the triangulation is finer and finer.
Motivated by the two disadvantages, we modify the classical discrete Gauss curvature $K_i$ to
\begin{equation}
R_i=\frac{K_i}{r_i^2}
\end{equation}
at the vertex $i$, which is called as $R$-curvature in the following.
In some sense, this new definition is more natural and respects the scaling exactly the same way as the Gauss curvature does.
Furthermore, this new definition can be used directly to approximate smooth Gauss curvature (see Section \ref{preliminaries and defnition of curvature}).
It is natural to consider the following combinatorial Yamabe problem for $R$-curvature.
\begin{question}\label{Question-2d-R-curvature}
Given a weighted triangulated surface $(M, \mathcal{T}, \Phi)$,
does there exist a circle packing metric $r$ that determines a constant $R$-curvature? How to find it if it exists?
\end{question}

We introduce the corresponding combinatorial versions of the Ricci flow and the Calabi flow for surfaces
and the Yamabe flow for 3-dimensional manifolds using the new definition of discrete curvature.
We give an answer to Question \ref{Question-2d-R-curvature} by discrete curvature flow methods and variational methods.
We just state the main results in this subsection.

We introduce a new combinatorial Ricci flow
\begin{equation}\label{equation of normalized combinatorial Ricci flow on surface introduction}
\frac{dg_i}{dt}=(R_{av}-R_i)g_i,
\end{equation}
where $R_{av}=2\pi\chi(M)/\|r\|^2_2$ and $g_i=r_i^2$ is a discrete version of Riemann metric tensor.
This flow has the same form as combinatorial Yamabe flow in two dimension.
By introducing a new discrete Laplace operator
\begin{equation*}\label{definition of Laplacian introduction}
\Delta f_i=\frac{1}{r_i^2}\sum_{j\sim i}(-\frac{\partial K_i}{\partial u_j})(f_j-f_i),
\end{equation*}
where $f$ is a function defined on all vertices and $u_j=\ln g_j$ is a coordinate transformation, we find that the flow (\ref{equation of normalized combinatorial Ricci flow on surface introduction}) exhibits similar properties to the smooth Ricci flow on surfaces. For example, the curvature $R_i$ evolves according to a heat-type equation
\begin{equation*}
\frac{dR_i}{dt}=\Delta R_i+R_i(R_i-R_{av}),
\end{equation*}
which has almost the same form as the evolution of the
scalar curvature along the classical Ricci flow and permits a discrete maximum principle.

We can also introduce a new combinatorial Calabi flow
\begin{equation}\label{equation of CCF introduction}
\frac{dg_i}{dt}=\Delta R_i\cdot g_i,
\end{equation}
which is similar to the smooth Calabi flow on surfaces.
By considering a discrete quadratic energy functional $\widetilde{C}(r)=\sum_{i=1}^{N}\varphi_{i}^2$, where $\varphi_{i}=K_i-R_{av}r_i^{2}$, we can define a modified Calabi flow
\begin{equation}\label{modified alpha Calabi flow introduction}
\dot{u}=-\frac{1}{2}\nabla_u \widetilde{C}.
\end{equation}

Combining Corollary \ref{global convergence of NCRF}, Theorem \ref{global convergence of CCF}, Theorem \ref{global convergence of modified CCF}, Theorem \ref{uniqueness of CCM} and Theorem \ref{main existence thm}, we have
the following main theorem.
\begin{theorem}\label{main convergence theorem in introduction}
Given a weighted triangulated surface $(M, \mathcal{T}, \Phi)$ with $\chi(M)\leq 0$.
Then the constant $R$-curvature metric is unique (if exists) up to scaling, and the following statements are mutually equivalent:
\begin{itemize}
  \item [(1)] There exists a constant $R$-curvature circle packing metric;
  \item [(2)] The solution to the combinatorial Ricci flow (\ref{equation of normalized combinatorial Ricci flow on surface introduction}) converges;
  \item [(3)] The solution to the combinatorial Calabi flow (\ref{equation of CCF introduction}) converges;
  \item [(4)] The solution to the modified Calabi flow (\ref{modified alpha Calabi flow introduction}) converges;
  \item [(5)] There exists a circle packing metric $r^*$ such that, for any nonempty proper subset $I$ of $V$,
\begin{equation}\label{combtopo-condition of GX-2-introduction}
2\pi\chi(M)\frac{\sum_{i\in I}r_i^{*2}}{\|r^*\|^{2}_{2}}>-\sum_{(e,v)\in Lk(I)}(\pi-\Phi(e))+2\pi\chi(F_I);
\end{equation}
  \item [(6)] The combinatorial-topological conditions $\mathscr{Y}\cap\mathds{R}^N_{<0}\neq\phi$ when $\chi(M)<0$, and $\mathscr{Y}\cap\{0\}\neq\phi$ when $\chi(M)=0$ are valid, where $\mathscr{Y}$ is defined in (\ref{Def-Y}).
\end{itemize}
\end{theorem}

\begin{remark}\label{Remark-two-example}
When $\chi(M)>0$, the conclusions in Theorem \ref{main convergence theorem in introduction} are not true.
Example \ref{example triangluate sphere} shows that the combinatorial Ricci flow for the sphere with the tetrahedron triangulation may not converge,
while Example \ref{example-not unique cccpm} shows that there are more than one constant curvature metric on the sphere with the tetrahedron triangulation.
\end{remark}

Theorem \ref{main convergence theorem in introduction} and Remark \ref{Remark-two-example}
provide a suitable answer to the Question \ref{Question-2d-R-curvature}.

As pointed out in Remark \ref{Remark-two-example},
it's surprising that the combinatorial flows are quite different from the smooth flows.
The combinatorial structure of the triangulation bring some extra trouble that needs special considerations.
However, for the surfaces with $\chi(M)>0$,
it's very interesting that we could evolve the flow (\ref{equation of normalized combinatorial Ricci flow on surface introduction})
very well by substituting the ``area element" $r_i^2$ by $r_i^{-2}$ or, more generally, $r_i^{\alpha}$ with $\alpha\leq0$.
In fact, for any $\alpha\in \mathds{R} $, we can define the $\alpha$-curvature as
\begin{equation}
R_{\alpha,i}=\frac{K_i}{r_i^{\alpha}}
\end{equation}
at each vertex $i$. The idea behind this definition is to consider $r^{\alpha}$ as a metric (of $\alpha$ order). From the viewpoint of Riemannian geometry, a piecewise flat metric is a singular Riemannian metric on $M^{2}$ or $M^{3}$, which produces conical singularities at all vertices. For any $\alpha\in\mathbb{R}$, a metric $g$ with
conical singularity at a point can be expressed as $g(z)=e^{f(z)}|z|^{2(\alpha-1)}dzd\bar{z}$ locally. Choosing $f(z)=-\ln\alpha^{2}$, then $g(z)=|dz^{\alpha}|^{2}$.
Comparing $r^{\alpha}$ with $|dz^{\alpha}|$, the $\alpha$-metric $r^{\alpha}$ may be taken as a discrete analogue of conical metric to some extent.
By introducing the corresponding $\alpha$-flows, we get the following result for constant $\alpha$-curvature problem, which is a combination of Theorem \ref{main alpha convergence theorem}, Theorem \ref{uniqueness for alpha curvature} and Corollary \ref{main existence thm for alpha-curvature}.

\begin{theorem}\label{main alpha convergence theorem in introduction}
Given a weighted triangulated surface $(M, \mathcal{T}, \Phi)$ with $\alpha\chi(M)\leq 0$. Then the constant $\alpha$-curvature metric is unique (if exists) up to scaling, and the following statements are mutually equivalent:
\begin{itemize}
  \item [(1)] There exists a constant $\alpha$-curvature circle packing metric;
  \item [(2)] The solution to the $\alpha$-Ricci flow (\ref{normalized alpha flow}) converges;
  \item [(3)] The solution to the $\alpha$-Calabi flow (\ref{definition of alpha Calabi flow}) converges;
  \item [(4)] The solution to the modified $\alpha$-Calabi flow (\ref{modified alpha Calabi flow}) converges;
  \item [(5)] There exists a circle packing metric $r^*$ such that, for any nonempty proper subset $I$ of $V$,
\begin{equation}\label{combtopo-condition of GX-2-introduction}
2\pi\chi(M)\frac{\sum_{i\in I}r_i^{*\alpha}}{\|r^*\|^{\alpha}_{\alpha}}>-\sum_{(e,v)\in Lk(I)}(\pi-\Phi(e))+2\pi\chi(F_I);
\end{equation}
  \item [(6)] The following combinatorial-topological conditions
              \begin{itemize}
                \item [(i)] $\mathscr{Y}\cap\mathds{R}^N_{<0}\neq\phi$, when $\alpha>0$ and $\chi(M)<0$;
                \item [(ii)] $\mathscr{Y}\cap\mathds{R}^N_{>0}\neq\phi$, when $\alpha<0$ and $\chi(M)>0$;
                \item [(iii)] $\frac{2\pi\chi(M)}{|V|}(1,\cdots,1)^T\in\mathscr{Y}$, when $\alpha\chi(M)=0$.
              \end{itemize}
              are valid, where $\mathscr{Y}$ is defined in (\ref{Def-Y}).
\end{itemize}
\end{theorem}

\begin{remark}
Theorem \ref{main alpha convergence theorem in introduction} contains Theorem \ref{Thm-thurston}, \ref{Thm-Chow-Luo}, \ref{Thm-Ge-Calabi}, \ref{main convergence theorem in introduction} as special cases.
Take $\alpha=2$, then Theorem \ref{main alpha convergence theorem in introduction} is reduced to Theorem \ref{main convergence theorem in introduction};
Take $\alpha=0$, then the equivalence between (1) and (5) is in fact Theorem \ref{Thm-thurston},
the equivalence between (1), (2) and (4) is in fact Theorem \ref{Thm-Chow-Luo} and Theorem \ref{Thm-Ge-Calabi}.
\end{remark}

Theorem \ref{main alpha convergence theorem in introduction} can not be generalized to
the case of $\alpha\chi(M)>0$ due to Remark \ref{Remark-two-example}.
so it is optimal to some extent.
We expect there are triangulations in this case that admits at least one constant $\alpha$-curvature metric, or more ideally, admits a solution to discrete Ricci flow that converges to a constant $\alpha$-curvature metric.
By the discrete maximum principle, we have the following existence result for non-negative constant $\alpha$-curvature metric.
\begin{theorem}\label{main theorem 4}
Suppose $(M, \mathcal{T}, \Phi)$ is a weighted triangulated surface. If there is a metric $r$ satisfying $R_{\alpha,i}\geq0$ for all $i\in V$, and
 \begin{equation*}
   -\sum_{(e,v)\in Lk(I)}(\pi-\Phi(e))+2\pi\chi(F_I)<0, \ \ \forall I:\emptyset \subsetneqq I \subsetneqq V,
 \end{equation*}
then there exists a non-negative constant $\alpha$-curvature metric $r^*$.
\end{theorem}

For a compact manifold $M^{3}$ with a triangulation $\mathcal{T}=\{V,E,F,T\}$, where $T$ is the set of all tetrahedrons.
A ball packing metric $r:V\rightarrow(0,+\infty)$ evaluates each $\{i,j\}\in E$ a length $l_{ij}=r_{i}+r_{j}$.
Denote $\alpha_{ijkl}$ as the solid angle at a vertex $i$ in the tetrahedron $\{i,j,k,l\}\in T$.
Cooper and Rivin \cite{CR} defined a discrete scalar curvature (we call it ``CR-curvature'' for simple in the following ) at each vertex $i$ as $$K_{i}=4\pi-\sum_{\{i,j,k,l\}\in T}\alpha_{ijkl}.$$
To study the constant CR-curvature problem, Glickenstein \cite{G1,G2} first introduced a combinatorial Yamabe flow
\begin{equation}
\frac{dr_{i}}{dt}=-K_{i}r_{i}.\label{glicken's-flow}
\end{equation}
Inspired by Glickenstein's work, the first author of this paper and Jiang \cite{GJ} modify Glickenstein's flow (\ref{glicken's-flow}) to
\begin{equation}
\frac{dr_{i}}{dt}=(\lambda-K_{i})r_{i},\label{Ge-jiang-flow}
\end{equation}
where $\lambda=\sum_{i=1}^N K_i r_i/\|r\|_{l^1}$, and get better convergence results.
Similar to the two dimensional case, we give a new definition of the combinatorial scalar curvature
$$R_i=\frac{K_i}{r_i^2}.$$
Similar to the smooth case, we find that the sphere packing metrics with constant combinatorial scalar curvature
are isolated and are exactly the critical points of a normalized Einstein-Hilbert-Regge functional $\sum R_ir_i^3/(\sum r_i^3)^{1/3}$.
Hence we propose to study a combinatorial Yamabe problem with respect to the new combinatorial scalar curvature $R_i$.
\begin{question}
Given a 3-dimensional manifold $M$ with triangulation $\mathcal{T}$, find a sphere packing metric with
constant combinatorial scalar curvature in the combinatorial
conformal class $\mathfrak{M}_{\mathcal{T}}$. Here $\mathfrak{M}_{\mathcal{T}}$ is the space of admissible sphere packing metrics determined by $\mathcal{T}$.
\end{question}

Set $R_{av}=\sum R_ir_i^3/\sum r_i^3$, we introduce the following combinatorial Yamabe flow
\begin{equation}\label{normalized comb Yamabe flow introduction}
\frac{dg_i}{dt}=(R_{av}-R_i)g_i
\end{equation}
to study the combinatorial Yamabe problem.
Nonpositive constant curvature metrics are local attractors of the flow (\ref{normalized comb Yamabe flow introduction}),
which implies the following result.
\begin{theorem}\label{convergence of CYF under existence introduction}
Suppose $r^*$ is a sphere packing metric on $(M, \mathcal{T})$ with nonpositive constant combinatorial
scalar curvature. If $||r(0)-r^*||^2$ is small enough, the solution of the normalized combinatorial
Yamabe flow (\ref{normalized comb Yamabe flow}) exists for all time and converges to $r^*$.
\end{theorem}

The paper is organized as follows. In Section \ref{preliminaries and defnition of curvature},
We introduce the new definition of discrete Gauss curvature for triangulated surfaces with
circle packing metrics.
In Section \ref{2-d comb RF}, We introduce the combinatorial Ricci flow and use it to
give a solution to the 2-dimensional Yamabe problem and
the prescribing curvature problem.
In Section \ref{Calabi flow section}, we introduce the combinatorial Calabi flow and study its properties.
In Section \ref{section-alpha}, we introduce the notion of $\alpha$-curvature and the corresponding $\alpha$-flows
and then study the corresponding constant curvature problem and prescribing curvature problem using the $\alpha$-flows.
In Section \ref{3-dimensional combinatorial Yamabe problem}, we introduce the new definition
of combinatorial scalar curvature on 3-dimensional manifolds and give a proof of Theorem \ref{convergence of CYF under existence introduction}.
In Section \ref{unsolved problems}, We list some unsolved problems closely related to the paper.

\section{The definition of combinatorial Gauss curvature}\label{preliminaries and defnition of curvature}\label{preliminaries of 2 dim}
The well-known discrete Gauss curvature $K_i$ has been widely studied in discrete geometry, we refer to \cite{CL1, GLSW, GGLSW} for recent progress.
However, there are two disadvantages of the classical definition of $K_i$.
For one thing, classical discrete Gauss curvature does not perform so perfectly in that it is scaling invariant, i.e. if $\tilde{r}_i=\lambda r_i$ for some positive constant $\lambda$, then $\tilde{K}_i=K_i$, which is different from the transformation of scalar curvature $R_{\lambda g}=\lambda ^{-1}R_g$ in the smooth case.
For another, classical discrete Gauss curvature can not be used directly to approximate smooth Gauss curvature when the triangulation is finer and finer.
As the triangulation of a fixed surface is finer and finer, we can get a sequence of polyhedrons which approximate the surface.
However, the classical discrete Gauss curvature at each vertex $i$ tends to zero,
for all triangles surrounding the vertex $i$ will finally run into the tangent plane
and the triangulation of the fixed surface locally becomes a triangulation of the tangent plane at vertex $i$.
To approximate the smooth Gauss curvature, dividing $K_i$ by an ``area elment" $A_i$ is necessary.
Since we are considering the circle packing metrics,
we may choose the area of the disk packed at $i$, i.e. $A_i=\pi r_i^2$, as the ``area elment" attached to the vertex $i$.
Omitting the coefficient $\pi$, we introduce the following definition of discrete curvature on triangulated surfaces.

\begin{definition}\label{defninition of curvature}
Given a weighted triangulated surface $(M, \mathcal{T}, \Phi)$ with circle packing metric $r: V\rightarrow (0,+\infty)$,
the discrete Gauss curvature at the vertex $i$ is defined to be
\begin{equation*}
R_i=\frac{K_i}{r_i^2},
\end{equation*}
where $K_i$ is the classical discrete Gauss curvature defined as the angle deficit at $i$ by (\ref{classical Gauss curv}).
\end{definition}

The curvature $R_i$ still locally measures the difference of the triangulated surface from the Euclidean plane at the vertex $i$.
In the rest of this subsection, we will give some evidences that
the Definition \ref{defninition of curvature} is a good candidate
to avoid the two disadvantages we mentioned in the above paragraph. Firstly, let us have a look at the following example.

\begin{figure}
\center
\begin{minipage}[t]{0.5\linewidth}
\centering
\includegraphics[width=0.4\textwidth]{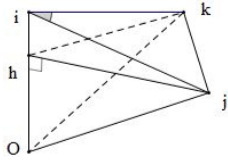}
\caption{{a local triangulation of sphere}}
\label{a local triangulation of sphere}
\end{minipage}
\end{figure}

\begin{example}
Consider the standard sphere $\mathbb{S}^2$ imbedded in $\mathds{R}^3$.
The smooth Gauss curvature of $\mathbb{S}^2$ is $+1$ everywhere.
Consider a local triangulation of $\mathbb{S}^2$ at the north pole, which is denoted as a vertex $i$.
$O$ is the origin of $\mathds{R}^3$. $h$ is a point lying in $Oi$, see Figure \ref{a local triangulation of sphere}.
The intersection between the horizonal plane passing through $h$ and $\mathbb{S}^2$ is a circle.
Divide this circle into $n$ equal portions.
Denote any two conjoint points by $j$ and $k$, thus $i$, $j$, $k$ forms a triangle.
Denote $x=\angle iOj$.  Taking $r_i$=$r_j$=$r_k$=$\overline{ij}/2$=$\sin\frac{x}{2}$, $\Phi_{ij}=\Phi_{ik}=0$
and suitably choosing $\Phi_{jk}$, we can always get a circle packing metric locally with $l_{jk}=2\sin \frac{\pi}{n}\sin x$.
Moreover, we can get an infinitesimal triangulation at vertex $i$ by letting $n\rightarrow +\infty$ and $x\rightarrow 0$.
By a trivial calculation we can get $\theta_i^{jk}=\angle kij=2\arcsin(sin\frac{\pi}{n}cos\frac{x}{2})$,
hence the classical discrete Gauss curvature is $K_i=2\pi-2n\arcsin(sin\frac{\pi}{n}cos\frac{x}{2})$. It's obviously that
$$\lim\limits_{n\rightarrow \infty,\,x\rightarrow 0}K_i=0.$$
However, the discrete Gauss curvature $R_i$ approaches the smooth Gauss curvature in the following way
$$\frac{R_i}{\pi}=\frac{K_i}{\pi r_i^2}=\frac{2\pi-2n\arcsin(sin\frac{\pi}{n}cos\frac{x}{2})}{\pi \sin^2\frac{x}{2}}\rightarrow +1.$$
\qed
\end{example}

Secondly, let's analyse from the view point of the scaling law of Riemann tensor and Gauss curvature. Recall that, for a $C^1$ curve $\gamma:[a, b]\rightarrow M$ in a Riemannian manifold $(M, g)$, the length of the curve is defined to be $L(\gamma, g)=\int_a^b\sqrt{g(\dot\gamma(t),\dot\gamma(t))}dt$, where $\dot\gamma(t)$ is the tangent vector of $\gamma(t)$ in $M$. So we have $L(\gamma, \tilde{g})=\lambda^{1/2}L(\gamma, g)$ if $\tilde{g}=\lambda g$ for some positive constant $\lambda$. Note that, for a weighted triangulated surface $(M, \mathcal{T}, \Phi)$ with circle packing metric $r: V\rightarrow (0,+\infty)$,
the length $l_{ij}$ of the edge $\{ij\}$ is given by (\ref{definition of length of edge}).
So the quantity corresponding to the Riemannian metric $g$
on the weighted triangulated surface $(M, \mathcal{T}, \Phi)$ with circle packing metric $r$
should be quadratic in $r$, among which $r_i^2$ is the simplest.
Furthermore, according to the definition $R_i=\frac{K_i}{r_i^2}$, we have $R_i(\lambda r_i^2)=\lambda^{-1}R_i(r_i^2)$, which has the same form as the transformation for the smooth
scalar curvature.

Finally, the Definition \ref{defninition of curvature} is geometrically reasonable. Let us recall the original definition of Gauss curvature.
We just give a sketch of the definition here and the readers could refer to \cite{S} for details.
Suppose $M$ is a surface embedded in $\mathds{R}^3$ and $\nu$ is the unit normal vector of $M$ in $\mathds{R}^3$.
$\nu$ defines the well-known Gauss map $\nu: M\rightarrow S^2\subset \mathds{R}^3$. Then the Gauss curvature at $p\in M$
is defined as
$$K(p)=\lim_{A\rightarrow p}\frac{Aera ~\nu(A)}{Aera ~ A},$$ where the limit is taken as the region $A$ around $p$ becomes smaller and smaller.
For the weighted triangulated surface $(M, \mathcal{T}, \Phi)$ with circle packing metric $r$, assume it is embedded in $\mathds{R}^3$ and take the normal vector as a set-valued function, then the definition of discrete Gauss curvature $R_i=\frac{K_i}{r_i^2}$ is an approximation of the original Gauss
curvature up to a uniform constant $\pi$, as the numerator $K_i$ is the measure of the set-valued function $\nu$ and the denominator $r_i^2$ is the area of the disk with radius $r_i$ up to a uniform constant $\pi$.

Note that the classical discrete Gauss curvature satisfies the Gauss-Bonnet identity (\ref{Gauss-Bonnet without weight}). If we define a discrete measure $\mu$ on the vertices by $\mu_i=r_i^2$, then we have
\begin{equation}\label{Def-discrete-measure}
\int_{M}fd\mu=\sum_{i=1}^N f_ir_i^2
\end{equation}
for any $f\in C(V)$. Using this discrete measure, we have the following discrete version of Gauss-Bonnet formula
\begin{equation}\label{Gauss-Bonnet with weight}
\int_MRd\mu=\sum_iR_id\mu_i=2\pi\chi(M).
\end{equation}
In this sense, the average curvature is
\begin{equation}\label{avarage curvature}
R_{av}=\frac{\int_MRd\mu}{\int_Md\mu}=\frac{2\pi\chi(M)}{||r||^2},
\end{equation}
where $||r||^2=\sum_{i=1}^Nr_i^2$ is the total measure of $M$ with respect to $\mu$.

We give the following table for the corresponding quantities for smooth surfaces and
triangulated surfaces.

\begin{table}[htbp]\label{table 1}
\centering  
\begin{tabular}{lccc}  
\hline
 &Smooth surfaces & Weighted triangulated surfaces\\ \hline  
Metric &$g=g_{ij}dx^idx^j$ & $g_i=r_i^2$ \\         
Length &$L(\gamma, g)=\int_a^b \sqrt{g(\dot\gamma(t),\dot\gamma(t))}dt$ & $l_{ij}=\sqrt{r_i^2+r_j^2+2r_ir_j\cos \Phi_{ij}}$ \\        
Measure &$d\mu=\sqrt{\det g_{ij}}dx$ & $d\mu_i=r_i^2$ \\
Curvature & Gauss curvature $K$ & $R_i=\frac{K_i}{r_i^2}$ \\
Total curvature & $\int_M Kd\mu=2\pi\chi(M) $ & $\sum_{i}R_id\mu_i=2\pi\chi(M)$ \\ \hline
\end{tabular}
\caption{Corresponding between smooth surfaces and triangulated surfaces}
\end{table}

For the curvature $R_i$ defined by (\ref{defninition of curvature}), it is natural to consider the
corresponding combinatorial Yamabe problem, i.e. does there exist any circle packing metric $r$ with constant
curvature $R_i$ ? Inspired by \cite{CL1, Ge, G1, L1},
we study the combinatorial Yamabe problem by the combinatorial curvature flows introduced in the following.

\section{The 2-dimensional combinatorial Ricci flow}\label{2-d comb RF}

\subsection{Definition of 2-dimensional combinatorial Ricci flow}

Ricci flow was introduced by Hamilton \cite{Ham2} to study the low dimension topology. For a closed Riemannian
manifold $(M^n, g_{ij})$, the Ricci flow
is defined as
\begin{equation}\label{Ricci flow equation}
\frac{\partial}{\partial t}g_{ij}=-2R_{ij}
\end{equation}
with normalization
\begin{equation}\label{normalized Ricci flow equation}
\frac{\partial}{\partial t}g_{ij}=\frac{2}{n}rg_{ij}-2R_{ij},
\end{equation}
where $r$ is the average of the scalar curvature. Ricci flow is a powerful tool
and has lots of applications. Specially, it has been used to prove
the Poincar\'{e} conjecture \cite{P1,P2,P3}. For a closed surface $(M^2, g_{ij})$, the Ricci flow
equation (\ref{Ricci flow equation}) is reduced to
\begin{equation}\label{Ricci flow equation on surface}
\frac{\partial}{\partial t}g_{ij}=-Rg_{ij}
\end{equation}
with the normalization
\begin{equation}\label{normalized Ricci flow equation on surface}
\frac{\partial}{\partial t}g_{ij}=(r-R)g_{ij},
\end{equation}
where $R$ is the scalar curvature. It is proved \cite{Ham1, CH1} that, for any closed
surface with any initial Riemannian metric, the solution of the normalized Ricci
flow (\ref{normalized Ricci flow equation on surface}) exists for all time
and converges to a constant curvature metric conformal to the initial metric as time
goes to infinity. It is further pointed out by Chen, Lu and Tian \cite{CLT}
that the Ricci flow can be used to give a
new proof of the famous uniformization theorem of Riemann surfaces. The discrete
version of surface Ricci flow was first introduce by Chow and Luo in \cite{CL1},
in which they gave another proof of Thurston's existence of circle packing theorem.
In this subsection, we will introduce a different combinatorial Ricci flow on surfaces
to study the constant curvature problem of $R_i$.

\begin{definition}\label{definition of modified combinatorial Ricci flow}
For a weighted triangulated surface $(M, \mathcal{T}, \Phi)$ with circle packing metric $r$,
the combinatorial Ricci flow is defined as
\begin{equation}\label{equation of modified combinatorial Ricci flow}
\frac{dg_i}{dt}=-R_ig_i,
\end{equation}
where $g_i=r_i^2$.
\end{definition}

Following Hamilton's approach, we introduce the following normalization of the flow (\ref{equation of modified combinatorial Ricci flow})
\begin{equation}\label{equation of normalized combinatorial Ricci flow on surface}
\frac{dg_i}{dt}=(R_{av}-R_i)g_i,
\end{equation}
where $R_{av}$ is the average curvature of $R_i$ defined by (\ref{avarage curvature}).
It is easy to check that the total measure $\mu(M)=||r||^2$ of $M$ with respect to $\mu$
is invariant along the normalized
flow (\ref{equation of normalized combinatorial Ricci flow on surface}), from which we know that
the average curvature $R_{av}=\frac{2\pi\chi(M)}{||r||^2}$ is invariant along
(\ref{equation of normalized combinatorial Ricci flow on surface}).
As our goal is to study the existence of constant combinatorial curvature metric, we
will focus on the properties of the flow (\ref{equation of normalized combinatorial Ricci flow on surface})
in the following. We assume $r(0)\in \mathbb{S}^{N-1}$ and then $r(t)\in \mathbb{S}^{N-1}$ along the flow in the following.

The flows (\ref{equation of modified combinatorial Ricci flow})
and (\ref{equation of normalized combinatorial Ricci flow on surface})
differ only by a change of scale in space and a change of parametrization in time.
Let $t, r, R$ denote the variables for the unnormalized flow (\ref{equation of modified combinatorial Ricci flow}),
and $\tilde{t}, \tilde{r}, \tilde{R}$ for the normalized flow (\ref{equation of normalized combinatorial Ricci flow on surface}).
Suppose $r(t)$, $t\in [0,T)$, is a solution of (\ref{equation of modified combinatorial Ricci flow}).
Set $\tilde{r}(\tilde{t})=\varphi(t)^{1/2}r(t)$, where $\varphi(t)=||r||^{-2}$ and $\tilde{t}=\int_0^t\varphi(\tau)d\tau$.
Then we have
$$||\tilde{r}||^2=1, \ \tilde{R}_i=\varphi(t)^{-1}R_i, \ \tilde{R}_{av}=\varphi(t)^{-1}R_{av}=2\pi\chi(M).$$
This gives
\begin{equation*}
\begin{aligned}
\frac{d\tilde{g}_i}{d\tilde{t}}=\frac{d\tilde{r}_i^2}{dt}\frac{dt}{d\tilde{t}}
=(\varphi'r_i^2+\varphi\frac{dr_i^2}{dt})\varphi^{-1}
=(\varphi'\varphi^{-1}\tilde{r}_i^2-\varphi R_i r_i^2)\varphi^{-1}
=(\tilde{R}_{av}-\tilde{R}_i)\tilde{g}_i,
\end{aligned}
\end{equation*}
where, in the last step, we use
$$\varphi'=\frac{d}{dt}||r||^{-2}=-||r||^{-4}\frac{d}{dt}||r||^2=2\pi\chi(M)\varphi^{2}.$$
Conversely, if $\tilde{r}(\tilde{t}), \tilde{t}\in [0, \tilde{T})$, is a solution of
(\ref{equation of normalized combinatorial Ricci flow on surface}),
set $r(t)=\varphi(\tilde{t})^{1/2}\tilde{r}(\tilde{t})$, where
$$\varphi(\tilde{t})=e^{-\tilde{R}_{av}\tilde{t}}, \ t=\int_0^{\tilde{t}}\varphi(\tau)d\tau.$$
Then it is easy to check that $\frac{d g_i}{dt}=-R_ig_i$.

For simplicity, we will set $u_i=\ln g_i$ in the following. Then the flows
(\ref{equation of modified combinatorial Ricci flow}) and (\ref{equation of normalized combinatorial Ricci flow on surface})
could be written as
\begin{equation}\label{equation of CRF in u}
\dot{u}=-R
\end{equation}
and
\begin{equation}\label{equation of NCRF in u}
\dot{u}=R_{av}\textbf{1}-R
\end{equation}
respectively, where $u=(u_1, \cdots, u_N)^T$, $\textbf{1}=(1, \cdots, 1)^T$ and $R=(R_1, \cdots, R_N)^T$.
\subsection{Discrete Laplace operator}

Since $\mu$ is an analogy of the area element, we can define an inner product $\langle\cdot, \cdot\rangle$ on
$(M, \mathcal{T}, \Phi)$ with circle packing metric $r$ by
\begin{equation}\label{inner product}
\langle f, h \rangle=\sum_{i=1}^{N}f_ih_ir_i^2=h^T\Sigma f
\end{equation}
for any real functions $f, h\in C(V)$, where $$\Sigma \triangleq diag\{r_1^2,\cdots,r_N^2\}.$$
A combinatorial operator $S: C(V)\rightarrow C(V)$ is said to be
self-adjoint if
$$\langle Sf, h\rangle=\langle f, Sh\rangle$$
for any $f, h\in C(V)$.

The classical discrete Laplace operator \cite{CHU} $``\Delta"$ is often written in the following form
$$\Delta f_i=\sum_{i\sim j}\omega_{ij}(f_j-f_i),$$
where the weight $\omega_{ij}$ can be arbitrarily selected for different purpose. Bennett Chow and Feng Luo \cite{CL1} first gave a special weight which comes from the dual structure of circle patterns, while Glickenstein \cite{G1} gave a similar weight in three dimension. The curvature $K_i$ is determined by the circle packing metric $r$. Consider the curvature as a function of $u$, where $u_i=\ln g_i=2\ln r_i$ are coordinate transformations, then Bennett Chow and Feng Luo's discrete Laplacian $\Delta_{CL}$ can be interpreted \cite{Ge} as the Jacobian of curvature map, i.e. $\Delta_{CL}=-2L$, where
\begin{equation}
L=\frac{\partial(K_1, \cdots, K_N)}{\partial (u_1,\cdots,u_N)}.
\end{equation}
It is noticeable that symbol $L$ above is different from that in \cite{Ge} by a factor 2,
which comes from the coordinate transformations.

Using the discrete measure $\mu_i=r_i^2$, we can give a new definition of discrete Laplace operator, which is slightly different from that
in \cite{CL1, Ge}.
\definition
For a weighted triangulated surface $(M, \mathcal{T}, \Phi)$ with circle packing metric $r: V\rightarrow (0,+\infty)$,
the discrete Laplace operator $\Delta:C(V)\rightarrow C(V)$ is defined as
\begin{equation}\label{definition of Laplacian}
\Delta f_i=\frac{1}{r_i^2}\sum_{j\sim i}(-\frac{\partial K_i}{\partial u_j})(f_j-f_i)
\end{equation}
for $f\in C(V)$.

The following property of $L$ will be used frequently in the rest of the paper.
\begin{lemma}\label{property of L}
(\cite{CL1}, Proposition 3.9.)
$L$ is positive semi-definite with rank $N$-$1$ and kernel $\{t\textbf{1}| t\in \mathds{R}\}$.
Moreover, $\frac{\partial K_i}{\partial u_i}>0$, $\frac{\partial K_i}{\partial u_j}<0$ for $i\sim j$ and $\frac{\partial K_i}{\partial u_j}=0$ for others.
\end{lemma}
Using this lemma, we can write discrete Laplace operator as a matrix, that is,
$$\Delta=-\Sigma^{-1} L$$
with
$$\Delta f=-\Sigma^{-1} Lf.$$
Note that the term $-\frac{\partial K_i}{\partial u_j}>0$  for $i\sim j$ in the definition of the operator $\Delta$,
thus it could be taken as a weight on the edges. We sometimes denote this weight as $w_{ij}$
in the following if there is no confusion. This implies that $\Delta$ is a standard Laplacian operator
defined on the weighted triangulated surface $(M, \mathcal{T}, \Phi)$ with circle packing metric $r$
and measure $\mu$. And it is easy to get
$$\langle \Delta f, g\rangle=g^T\Sigma(-\Sigma^{-1} Lf)=-g^TLf=\langle f, \Delta g\rangle.$$
Hence the Laplacian operator $\Delta$ is self-adjoint with respect to $\langle\cdot, \cdot\rangle$. Moreover, if we define
$$|\nabla f|_i^2=\frac{1}{2r_i^2}\sum_{j\thicksim i}(-L_{ij})(f_j-f_i)^2,$$
then we can derive
\begin{equation}\label{integration by part}
\langle |\nabla f|^2, \textbf{1} \rangle=\sum_{i=1}^N|\nabla f|_i^2r_i^2=\sum_{j\thicksim i}(-L_{ij})(f_j-f_i)^2=f^TLf=-\langle f, \Delta f\rangle.
\end{equation}
Using the discrete measure $\mu$ introduced in (\ref{Def-discrete-measure}), (\ref{integration by part}) can be written in an integral form
$$\int_{M}|\nabla f|^2d\mu=-\int_{M} f\Delta fd\mu.$$

When the circle packing metric evolves along the flow (\ref{equation of normalized combinatorial Ricci flow on surface}),
so is the curvature $R_i$. The evolution of $R_i$ is very simple and it has almost the same
form as the evolution of scalar curvature along the Ricci flow on surfaces derived by Hamilton in \cite{Ham1}.

\begin{lemma}\label{lemma for evolution of R_i under NCRF}
Along the normalized combinatorial Ricci flow (\ref{equation of normalized combinatorial Ricci flow on surface}),
the curvature $R_i$ satisfies the evolution equation
\begin{equation}\label{evolution of R_i under NCRF}
\frac{dR_i}{dt}=\Delta R_i+R_i(R_i-R_{av}).
\end{equation}
\end{lemma}
\textbf{Proof.}
As $u_i=\ln r_i^2$, we have $\frac{\partial}{\partial u_j}=\frac{1}{2}r_j\frac{\partial}{\partial r_j}$. Then we have
\begin{equation*}
  \begin{aligned}
  \frac{dR_i}{dt}
  &=\sum_j\frac{\partial R_i}{\partial u_j}\frac{du_j}{dt}\\
  &=\sum_j(\frac{1}{r_i^2}\frac{\partial K_i}{\partial u_j}-\frac{K_i}{r_i^4}r_ir_j\delta_{ij})(R_{av}-R_j)\\
  &=\frac{1}{r_i^2}\sum_j\frac{\partial K_i}{\partial u_j}(R_{av}-R_j)-R_i(R_{av}-R_i)\\
  &=-\frac{1}{r_i^2}\sum_j\frac{\partial K_i}{\partial u_j}R_j+R_i(R_i-R_{av})\\
  &=\frac{1}{r_i^2}\sum_{j\sim i}(-\frac{\partial K_i}{\partial u_j})(R_j-R_i)+R_i(R_i-R_{av}).
  \end{aligned}
\end{equation*}
In the last two steps, Lemma \ref{property of L} is used. \qed\\

Thus the evolution equation (\ref{evolution of R_i under NCRF}) of $R_i$
is a reaction-diffusion equation. In fact, suppose $\Omega$ is a subset of $V$, set $\nabla_{ij}f=f_j-f_i$ for $f\in C(V)$
and $i\sim j$, and denote
$$\partial \Omega=\{\{ij\}\in E|i\in \Omega, j\in \Omega^c\},$$
then it is easy to check that
$$\frac{d}{dt}\int_{\Omega}Rd\mu=\sum_{i\in \Omega, j\in \Omega^c, i\sim j}(-\frac{\partial K_i}{\partial u_j})(R_j-R_i)
=\int_{\partial\Omega}\nabla R dw.$$
Then we have similar explanation as the smooth Ricci flow for equation (\ref{equation of normalized combinatorial Ricci flow on surface}).

\subsection{Maximum principle}
As the evolution equation (\ref{evolution of R_i under NCRF}) is a heat-type equation, we have the following
discrete maximum principle for such equations. Such maximum principle is almost standard, for completeness,
we write it down here.

\begin{theorem}\label{Maximum priciple}(Maximum Principle)
Let $f: V\times [0, T)\rightarrow \mathds{R}$ be a $C^1$ function such that
$$\frac{\partial f_i}{\partial t}\geq \Delta f_i+ \Phi_i(f_i), \ \ \forall (i, t)\in V\times [0, T) $$
where
$\Phi_i: \mathds{R}\rightarrow \mathds{R}$ is a local Lipschitz function.
Suppose there exists $C_1\in \mathds{R}$ such that $f_i(0)\geq C_1$ for all $i\in V$. Let $\varphi$ be the
solution to the associated ODE
\begin{equation*}
\begin{aligned}
\left\{
  \begin{array}{ll}
    \frac{d \varphi}{dt}=\Phi_i(\varphi) \\
    \varphi(0)=C_1,
  \end{array}
\right.
\end{aligned}
\end{equation*}
then
$$f_i(t)\geq \varphi(t)$$
for all $(i, t)\in V\times [0, T)$ such that $\varphi(t)$ exists.

Similarly, suppose $f: V\times [0, T)\rightarrow \mathds{R}$ be a $C^1$ function such that
$$\frac{\partial f_i}{\partial t}\leq \Delta f_i+ \Phi_i(f_i), \ \ \forall (i, t)\in V\times [0, T). $$
Suppose there exists $C_2\in \mathds{R}$ such that $f_i(0)\leq C_2$ for all $i\in V$. Let $\psi$ be the
solution to the associated ODE
\begin{equation*}
\begin{aligned}
\left\{
  \begin{array}{ll}
    \frac{d \psi}{dt}=\Phi_i(\psi) \\
    \psi(0)=C_2,
  \end{array}
\right.
\end{aligned}
\end{equation*}
then
$$f_i(t)\leq \psi(t) $$
for all $(i, t)\in V\times [0, T)$ such that $\psi(t)$ exists.
\end{theorem}

\begin{remark}
In fact, Theorem \ref{Maximum priciple} is valid for general Laplacian operators
defined as
$$\Delta f_i=\sum_{j\sim i}a_{ij}(t)(f_j-f_i),$$
where $a_{ij}\geq 0$, but not required to satisfy the symmetry condition $a_{ij}=a_{ji}$.
\end{remark}

The proof of Theorem \ref{Maximum priciple} is almost the same as that in \cite{CK}, we give a proof of it here just for completeness.

\begin{proposition}
Let $f: V\times [0,T)\rightarrow \mathds{R}^N$ be a $C^1$ solution to the heat equation
$$\frac{\partial f}{\partial t}=\Delta f,$$
i.e.
$$\frac{\partial f_i}{\partial t}=\Delta f_i=\sum_{j\sim i}a_{ij}(t)(f_j-f_i).$$
If there are constants $C_1\leq C_2\in \mathds{R}$ such that $C_1\leq f_i(0)\leq C_2$ for all $i\in V$, then
$$C_1\leq f_i(t)\leq C_2,\ \ \  \forall (i, t)\in V\times [0,T).$$
\end{proposition}
This proposition follows from the following more general theorem.

\begin{theorem}\label{MP for supersolution without reaction term}
Suppose $f:V\times [0, T)\rightarrow \mathds{R}^N$ be a $C^1$ function with $f_i(0)\geq \alpha$ for
some constant $\alpha\in \mathds{R}$, and $f$ satisfies
$$\frac{\partial f_i}{\partial t}\geq \Delta f_i$$
at any $(i, t)\in V\times [0,T)$ such that $f_i(t)<\alpha$. Then
$$f_i(t)\geq \alpha$$
for all $(i, t)\in V\times [0, T)$.
\end{theorem}

\proof
Note that if $H:V\times [0, T)\rightarrow \mathds{R}$ is a $C^1$ function and $(i, t_0)$ is the
vertex and time where $H_i(t)$ attains its minimum among all vertices and earlier times, namely
$$H_i(t_0)=\min _{(j, t)\in V\times [0, t_0]}H_j(t),$$
then
$$\frac{\partial H_i}{\partial t}(t_0)\leq 0,  \Delta H_i(t_0)\geq 0.$$
Consider a function $H$ defined by
$H_i(t)=\left(f_i(t)-\alpha\right)+\varepsilon t +\varepsilon,$
where $\varepsilon$ is any positive constant. Then $H_i(0)\geq \varepsilon>0$, and
$$\frac{\partial H_i}{\partial t}=\frac{\partial f_i}{\partial t}+\varepsilon\geq \Delta f_i+\varepsilon=\Delta H_i+\varepsilon$$
for $(i, t)\in V\times [0,T)$ satisfies $f_i(t)\leq \alpha$. Thus we just need to prove that $H_i(t)>0$ for $(i, t)\in V\times [0,T)$.

Suppose that $(i, t_0)$ is the first point in $V\times [0, T)$ such that $H_i(t_0)=0$. Then
$$f_i(t_0)=H_i(t_0)+\alpha-\varepsilon t_0-\varepsilon=\alpha-\varepsilon t_0-\varepsilon<\alpha.$$
Then we have
$$0\geq \frac{\partial H_i}{\partial t}(t_0)\geq \Delta H_i(t_0)+\varepsilon\geq \varepsilon>0,$$
which is a contradiction. And the theorem follows from $\varepsilon\rightarrow 0$. \qed

Now suppose $\beta: V\times [0, T)\rightarrow \mathds{R}$ be a function.

\begin{proposition}\label{MP with linear reaction term}
Let $f: V\times [0, T) \rightarrow \mathds{R}$ be a $C^1$ function such that
$$\frac{\partial f_i}{\partial t}\geq \Delta f_i+\beta_if_i.$$
And for each $\tau\in [0, T)$, $\exists C_\tau< \infty$ such that $\beta_i(t)<C_\tau$ for $(i, t)\in V\times [0, \tau]$.
If $f_i\geq 0$ for all $i\in V$, then $f_i(t)\geq 0$ for all $(i, t)\in V\times [0, T)$.
\end{proposition}
\proof
Given $\tau\in (0, T)$, define
$$J_i(t)=e^{-C_\tau t}f_i(t),$$
then
$$\frac{\partial J_i}{\partial t}\geq \Delta J_i+(\beta_i-C_\tau)J_i.$$
Since $\beta_i(t)\leq C_\tau$ for all $(i, t)\in V\times [0, \tau]$, then we have
$$\frac{\partial J_i}{\partial t}\geq \Delta J_i$$
for $(i, t)\in V\times [0, \tau]$ such that $J_i(t)<0$. Then by the Theorem \ref{MP for supersolution without reaction term},
we have $J_i(t)\geq 0$, and then $f_i(t)\geq 0$ for any $(i, t)\in V\times [0, T)$. \qed

Now we give the proof of Theorem \ref{Maximum priciple}.\\

\textbf{Proof of Theorem \ref{Maximum priciple}:}
We take the lower bound to prove, and proof for the upper bound is similar.
By assumption, we have
$$\frac{\partial}{\partial t}(f_i-\varphi)\geq \Delta (f_i-\varphi)+\Phi_i(f_i)-\Phi_i(\varphi).$$
The assumptions on the initial data imply that
$$(f_i-\varphi)(0)\geq 0, \forall i\in V.$$
Let $\tau\in (0, T)$, then $\exists C_\tau<\infty$ such that $|f_i(t)|<C_\tau$ for all $(i, t)\in V\times [0, \tau]$
and $|\varphi(t)|<C_\tau$ for $t\in [0, \tau]$. Since $F$ is locally Lipschitz, then we have
$$|\Phi_i(f_i)-\Phi_i(\varphi)|\leq L_\tau |f_i-\varphi|, \forall i\in V$$
for some positive constant $L_\tau>0$.
Then
$$\frac{\partial}{\partial t}(f_i-\varphi)\geq \Delta (f_i-\varphi)-L_\tau Sgn(f_i-\varphi)\cdot (f_i-\varphi),
\forall (i, t)\in V\times [0, \tau].$$
Applying Proposition \ref{MP with linear reaction term}, we get
$$f_i(t)-\varphi(t)\geq 0, \forall (i, t)\in V\times [0, \tau].$$
And the theorem follows from $\tau\in [0, T)$ is arbitrary. \qed

Applying Theorem \ref{Maximum priciple} to (\ref{evolution of R_i under NCRF}) with vanishing initial value, we can easily get the following corollary.

\begin{corollary}\label{preserve positive curv}
If $R_i(0)\geq 0$ ($R_i(0)\leq 0$) for all $i\in V$, then $R_i(t)\geq 0$ ($R_i(t)\leq 0$) for  all $i\in V$ as long as the flow exists,
i.e. the positive and negative curvatures are preserved along the normalized flow (\ref{equation of normalized combinatorial Ricci flow on surface}).
\end{corollary}

Set $R_{max}(t)=\max_{i\in V}R_i(t)$ and $R_{min}(t)=\min_{i\in V}R_i(t)$.
Applying Theorem \ref{Maximum priciple} to (\ref{evolution of R_i under NCRF})
with general initial data, we get the following lower bounds
for $R_i(t)$. For simplicity, We just state it in the form parallel to that given in \cite{CK}.

\begin{lemma}\label{lower bound}(Lower Bound)
Let $r_i(t)$ be a solution to the normalized combinatorial Ricci flow (\ref{evolution of R_i under NCRF})
on a closed triangulated surface $(M, T, \Phi)$.
\begin{itemize}
  \item[(1)] If $\chi(M)<0$, then
$$R_i-R_{av}\geq \frac{R_{av}}{1-(1-\frac{R_{av}}{R_{min}(0)})e^{R_{av}t}}-R_{av}\geq (R_{min}(0)-R_{av})e^{R_{av}t}.$$
  \item[(2)] If $\chi(M)=0$, then
$$R_i\geq \frac{R_{min}(0)}{1-R_{min}(0)t}>-\frac{1}{t}.$$
  \item[(3)] If $\chi(M)>0$ and $R_{min}(0)<0$, then
$$R\geq \frac{R_{av}}{1-(1-\frac{R_{av}}{R_{min}(0)})e^{R_{av}t}}\geq R_{min}(0)e^{-R_{av}t}.$$
\end{itemize}
\end{lemma}

Notice that in each case, the right hand side of the lower bound estimate
tends to 0 as $t\rightarrow \infty$. However, for the upper bound, the situation is not so good.
The solution of the ODE
\begin{equation*}
\begin{aligned}
\left\{
  \begin{array}{ll}
    \frac{ds}{dt}&=s(s-R_{av})  \\
    s(0)&=s_0
  \end{array}
\right.
\end{aligned}
\end{equation*}
corresponding to (\ref{evolution of R_i under NCRF}) is
\begin{equation*}
\begin{aligned}
s(t)=\left\{
  \begin{array}{ll}
    0, & \hbox{$s_0=0$} \\
    \frac{s_0}{1-s_0t}, & \hbox{$s_0\neq 0,\ \chi(M)=0$} \\
    \frac{R_{av}}{1-(1-\frac{R_{av}}{s_0})e^{R_{av}t}}, & \hbox{$s_0\neq 0,\ \chi(M)\neq 0$}
  \end{array}
\right.
\end{aligned}.
\end{equation*}
If $s_0>\max \{R_{av}, 0\}$, there is $T<\infty$ given by
\begin{equation*}
\begin{aligned}
T=\left\{
    \begin{array}{ll}
      -\frac{1}{R_{av}}\ln \left(1-\frac{R_{av}}{s_0}\right)>0, & \hbox{$\chi(M)\neq 0$} \\
      \frac{1}{s_0}, & \hbox{$\chi(M)=0$}
    \end{array}
  \right.
\end{aligned}
\end{equation*}
such that
\begin{equation*}
\lim_{t\rightarrow T^-}s(t)=+\infty.
\end{equation*}
The implies that, in the case of $R_{max}(0)>\max\{R_{av}, 0\}$, Theorem \ref{Maximum priciple}
will not give us good upper
bounds for the curvature along the combinatorial Ricci flow (\ref{evolution of R_i under NCRF}).
However, in the case of $R_{max}(0)<0$, we have good upper bounds.

\begin{lemma}\label{upper bound for R_{i}(0)<0}
If $R_{i}(0)<0$ for all $i\in V$, then we have
$$R_i(t)-R_{av}\leq R_{av}\left(1-\frac{R_{av}}{R_{max}(0)}\right)e^{R_{av}t}.$$
\end{lemma}

These results are almost parallel to that of the smooth Ricci flow on surfaces.

\subsection{Convergence of the 2-dimensional combinatorial Ricci flow}
Combining Lemma (\ref{lower bound}) with (\ref{upper bound for R_{i}(0)<0}) gives
\begin{theorem}\label{existence theorem for negative curvature}
Suppose $(M, \mathcal{T}, \Phi)$ is a weighted triangulated surface
with initial circle packing metric $r(0)$ satisfying $R_{i}(0)<0$ for all $i\in V$,
then there exists a negative constant curvature metric $r^*$ on $(M, \mathcal{T}, \Phi)$.
Furthermore, the solution $r(t)$ to the normalized combinatorial Ricci flow
(\ref{equation of normalized combinatorial Ricci flow on surface}) exists for all time and
converges exponentially fast to $r^*$ as $t\rightarrow +\infty$.
\end{theorem}

Theorem \ref{existence theorem for negative curvature}, derived by combinatorial maximum principle, implies more than it seems.
It claims that, if there is a metric with all curvatures negative,
then there always exists a negative constant curvature metric, and vice versa.
Using this fact, we will give a combinatorial and topological condition
which is equivalent to the existence of negative constant curvature metric in
subsection \ref{section-existence and uniqueness of cccpm}.
For the convergence of combinatorial Ricci flow,
the negative initial curvature condition in Theorem \ref{existence theorem for negative curvature} is a little restrictive.
It is natural to consider Hamilton's approach to generalize this result. However, we found that there is some technical difficulties to go further
with Hamilton's approach in this case. By studying critical points of (\ref{equation of normalized combinatorial Ricci flow on surface}),
which is an ODE system, we can get the following local convergence result.

\begin{theorem}\label{local convergence for NCRF under general conditions}
Suppose  $r^*$ is a constant $R$-curvature metric on a weighted triangulated surface $(M, \mathcal{T}, \Phi)$.
If the first positive eigenvalue $\lambda_1$ of $-\Delta$ at $r^*$ satisfies
\begin{equation}\label{lambda(lambda)>c general condition}
\lambda_1(-\Delta)>R_{av}^*=\frac{2\pi\chi(M)}{||r^*||^2}
\end{equation}
and $||r(0)-r^*||^2$ is small enough, then the solution to the normalized combinatorial
Ricci flow (\ref{equation of normalized combinatorial Ricci flow on surface}) exists for $t\in [0,+\infty)$ and converges exponentially fast to $r^*$.
\end{theorem}

\textbf{Proof.}
We can rewrite the normalized combinatorial Ricci flow (\ref{equation of normalized combinatorial Ricci flow on surface})
as
$$\frac{dr_i}{dt}=\frac{1}{2}(R_{av}-R_i)r_i.$$
Set $\Gamma_i(r)=\frac{1}{2}(R_{av}-R_i)r_i$, then the Jacobian matrix of $\Gamma(r)$ is given by
\begin{equation*}
(D_r\Gamma(r))_{ij}
=\frac{\partial}{\partial r_j}\left(\frac{1}{2}(R_{av}-R_i)r_i\right)=\frac{1}{2}(R_{av}-R_i)\delta_{ij}+
(R_i\delta_{ij}-R_{av}\frac{r_ir_j}{||r||^2})-\frac{1}{r_ir_j}\frac{\partial K_i}{\partial u_j}.
\end{equation*}
Denote $\Lambda=\Sigma^{-\frac{1}{2}} L \Sigma^{-\frac{1}{2}}$, then $\Lambda\sim \Sigma^{-\frac{1}{2}}\Lambda\Sigma^{\frac{1}{2}}=-\Delta$, hence $\lambda_1(-\Delta)=\lambda_1(\Lambda)$. At the constant curvature metric point $r^*$, we have
\begin{equation*}
D_r\Gamma|_{r^*}=R_{av}\left(I-\frac{r r^T}{||r||^2}\right)-\Sigma^{-\frac{1}{2}} L \Sigma^{-\frac{1}{2}}=R_{av}\left(I-\frac{r r^T}{||r||^2}\right)-\Lambda.
\end{equation*}
Select an orthonormal matrix $P$ such that
$$P^T\Lambda P=diag\{0,\lambda_1(\Lambda),\cdots,\lambda_{N-1}(\Lambda)\}.$$
Suppose $P=(e_0,e_1,\cdots,e_{N-1})$,
where $e_i$ is the $(i+1)$-column of $P$.
Then $\Lambda e_0=0$ and $\Lambda e_i=\lambda_i e_i,\,1\leq i\leq N-1$,
which implies $e_0=r/\|r\|$ and $r\perp e_i,\,1\leq i\leq N-1$.
Hence $\big(I_{N}-\frac{rr^T}{\|r\|^{2}}\big)e_0=0$ and $\big(I_{N}-\frac{rr^T}{\|r\|^{2}}\big)e_i=e_i$, $1\leq i\leq N-1$,
which implies $P^T\big(I_{N}-\frac{rr^T}{\|r\|^{2}}\big)P=diag\{0,1,\cdots,1\}$. Therefore,
\begin{equation*}
D_{r}\Gamma\big|_{r^*}=
P \cdot diag\{0,R_{av}-\lambda_1(\Lambda),\cdots,R_{av}-\lambda_{N-1}(\Lambda)\}\cdot P^T.
\end{equation*}
If $\lambda_1(\Lambda)>R_{av}^*=\frac{2\pi\chi(M)}{||r^*||^2}$, then $D_{r}\Gamma\big|_{r^*}$ is negative semi-definite with kernel $\{cr| c\in \mathbb{R}\}$ and $rank\,(D_{r}\Gamma\big|_{r^*})=N-1$. Note that, along the flow (\ref{equation of normalized combinatorial Ricci flow on surface}), $||r||^2$ is invariant. Thus the kernel is transversal to the flow. This implies that $D_r\Gamma|_{r^*}$ is negative definite on $\mathbb{S}^{N-1}$ and $r^*$ is a local attractor of the normalized combinatorial Ricci flow (\ref{equation of normalized combinatorial Ricci flow on surface}). Then the conclusion follows from the Lyapunov Stability Theorem(\cite{P}, Chapter 5).\qed\\

If $\chi(M)\leq 0$, we always have $\lambda_1(-\Delta)>0\geq R_{av}^*=\frac{2\pi\chi(M)}{||r^*||^2}$, thus we get
\begin{corollary}\label{local convergence for NCRF}
Suppose  there is a nonpositive constant curvature metric $r^*$ on a weighted triangulated surface $(M, \mathcal{T}, \Phi)$ with
$\chi(M)\leq 0$. If $||r(0)-r^*||^2$ is small enough, then
the solution to the normalized combinatorial Ricci flow (\ref{equation of normalized combinatorial Ricci flow on surface})
exists for $t\in [0,+\infty)$ and converges exponentially fast to $r^*$.
\end{corollary}

Before proving the global results, we first check that the convergence of
the flow (\ref{equation of normalized combinatorial Ricci flow on surface}) ensures the existence of the constant curvature metric.
More generally, we have the following result.

\begin{lemma}\label{compact imply CCCPM exist}
Suppose that the solution to the flow (\ref{equation of normalized combinatorial Ricci flow on surface})
lies in a compact region in $\mathds{R}^N_{>0}$, then there exists a constant curvature metric $r^*$.
Moreover, there exists a sequence of metrics which converge to $r^*$ along this flow.
\end{lemma}

\textbf{Proof.}
Consider the Ricci potential
\begin{equation*}
F(u)=\int_{u_0}^u\sum_{i=1}^{N}(K_i-R_{av}r_i^2)du_i,
\end{equation*}
where $u_i=\ln r_i^2$ and $u_0\in \mathds{R}^N$ is an arbitrary point. If we set $\varphi_i=K_i-\frac{2\pi \chi(M)}{||r||^2}r_i^2$, by direct calculations, we have $$\frac{\partial \varphi_i}{\partial u_j}=\frac{\partial \varphi_j}{\partial u_i}.$$
As the domain of $u$ is simply connected, this implies that the integration is path independent and then the Ricci potential is well defined.
Suppose $t\in[0,T)$ and $T$ is the maximal existing time of $u(t)$. Then $\{u(t)\}\subset\subset \mathds{R}^{N}$ implies $T=+\infty$, otherwise $u(t)$ will run out of the compact region. $\{u(t)\}\subset\subset \mathds{R}^{N}$ also implies that $F(u(t))$ is bounded. Moreover,
\begin{equation*}
\frac{d}{dt}F(u(t))=(\nabla_u F)^T\cdot \dot{u}=-\sum\left(\frac{K_i}{r_i}-\frac{2\pi \chi(M)}{||r||^2}r_i\right)^2=-\sum_i(R_i-R_{av})^2r_i^2\leq 0.
\end{equation*}
Hence $\lim_{t\rightarrow +\infty}F(u(t))$ exists.
Then by the mean value theorem, there exists a sequence $\xi_n\in (n, n+1)$ such that
$$F(u(n+1))-F(u(n))=(F(u(t)))'|_{t=\xi_n}
=-\sum_i(R_i-R_{av})^2r_i^2(\xi_n)\rightarrow 0$$
as $n\rightarrow +\infty$.
Since $\{r(t)\}\subset\subset R_{>0}^N$,  we can further select a subsequence $t_k=\xi_{n_k}$ such that $r(t_k)\rightarrow r^*$ as $k\rightarrow +\infty$. Combining $\sum_ir_i^2(R_i-R_{av})^2(\xi_{n_k})\rightarrow 0$, we know that $r^*$ determines a constant $R$-curvature.
\qed

\begin{corollary}\label{convergence imply CCCPM exist}
Suppose that the solution to flow (\ref{equation of normalized combinatorial Ricci flow on surface}) exists for all time and converges to $r(+\infty)$, then $r(+\infty)$ is a constant curvature metric.
\end{corollary}

\begin{remark}
For the Ricci potential functional $F$ used in the proof of Lemma \ref{compact imply CCCPM exist},
the primitive form of this type of functional was first constructed by Colin de Verdi\`{e}re in \cite{DV} and then
further studied by Chow and Luo in \cite{CL1}. Our definition is a modification of their constructions.
\end{remark}

For the Ricci potential $F$ introduced in the proof of Lemma \ref{compact imply CCCPM exist}, we further have the following
property.

\begin{lemma}\label{property of Ricci potential under general condition}
Given a weighted triangulated surface $(M, \mathcal{T}, \Phi)$.
Assume there is a constant curvature circle packing metric $r^*$ on $(M, \mathcal{T}, \Phi)$. Define the Ricci potential as
\begin{equation}\label{Ricci potential}
F(u)=\int_{u^*}^u\sum_{i=1}^{N}(K_i-R_{av}r_i^2)du_i.
\end{equation}
Denote
$$\mathscr{U}_a\triangleq\{u\in \mathds{R}^N|\sum_i u_i=a\},\ \ \ a\in \mathds{R}.$$
If $\lambda_1(-\Delta)>R_{av}$ for all $r\in \mathds{R}_{>0}^N$, then
\begin{itemize}
  \item[(1)] $Hess_uF$ is positive semi-definite with rank $N-1$ and kernel $\{t \mathbf{1}|t\in \mathds{R}\}$.
  \item[(2)] Restricted to $\mathscr{U}_a$, $F|_{\mathscr{U}_a}$ is strictly convex and proper.
   $F|_{\mathscr{U}_a}$ has a unique zero point which is also the unique minimum point. Moreover, $\lim\limits_{u\in\mathscr{U}_a,\,u\rightarrow \infty}F(u)=+\infty$.
\end{itemize}
\end{lemma}
\textbf{Proof.}
By direct calculations, we have
\begin{equation}\label{Hession of F}
\begin{aligned}
Hess_uF
=L-R_{av}\Sigma^{\frac{1}{2}}
\left(I-\frac{rr^T}{||r||^2}\right)\Sigma^{\frac{1}{2}}=\Sigma^{\frac{1}{2}}
\left(\Lambda-R_{av}\left(I-\frac{rr^T}{||r||^2}\right)\right)
\Sigma^{\frac{1}{2}}.
\end{aligned}
\end{equation}
By the same analysis as that of Theorem \ref{local convergence for NCRF under general conditions},
we know that, if $\lambda_1(-\Delta)>R_{av}$, $Hess_u F$ is positive semi-definite with
$rank(Hess_u F)=N-1$ and $Ker(Hess_uF)=\{t \mathbf{1}|t\in \mathds{R}\}$.

For the second part of the proof, we follow that of Colin de Verdi\`{e}re \cite{DV} and  Chow and Luo\cite{CL1}.
It's easy to see $F(u)=F(u+t\textbf{1})$ for any $t\in \mathds{R}$ and $u\in \mathds{R}^N$.
As $F$ is invariant along the direction $\textbf{1}$, we just need to prove it on $\mathscr{U}_0=\{u\in\mathds{R}^N|\sum u_i=0\}$.
A rigorous proof is formulated in Theorem B.2 of \cite{Ge}.
In addition, $u^*+\frac{a-\textbf{1}^Tu^*}{N}\textbf{1}$ is the unique zero point and minimum point of $F|_{\mathscr{U}_a}$.\qed

\begin{theorem}\label{global convergence of NCRF under general condition}
Suppose $(M, \mathcal{T}, \Phi)$ is a weighted triangulated surface and $\lambda_1(-\Delta)>R_{av}$ for all $r\in \mathds{R}_{>0}^N$.
Then the solution to the flow (\ref{equation of normalized combinatorial Ricci flow on surface}) converges if and only if there exists a constant curvature metric $r^*$. Furthermore, if the solution converges, it converges exponentially fast to the metric of constant curvature.
\end{theorem}
\textbf{Proof.}
The necessary part can be seen from Corollary \ref{convergence imply CCCPM exist}.
For the sufficient part, without loss of generality, we may assume $r^*\in\mathbb{S}^{N-1}$.
Consider the flow (\ref{equation of normalized combinatorial Ricci flow on surface})
with initial metric $r(0)\in\mathbb{S}^{N-1}\cap\mathds{R}^N_{>0}$.
Then $r(t)\in\mathbb{S}^{N-1}\cap\mathds{R}^N_{>0}$ and $u(t)\in \mathscr{V}_1^*\triangleq\{u\in \mathds{R}^N|\sum_i e^{u_i}=1\}$ for all $t$.
Let $\Pi$ be the orthogonal projection from $\mathscr{V}_1^*$ to the plane $\mathscr{U}_0=\{u\in \mathds{R}^N|\sum_i u_i=0\}$.
Then $F(u)=F(\Pi(u))$ and $\Pi(u)\rightarrow \infty$ as $u\rightarrow \infty$.
In fact, for a sequence $\{u^{(n)}\}$, $\Pi(u^{(n)})$ is unbounded
if and only if $|\Pi(u^{(n)})_i-\Pi(u^{(n)})_j|=|u^{(n)}_i-u^{(n)}_j|$ is unbounded.
Then Lemma \ref{property of Ricci potential under general condition} implies
\begin{equation}
\lim\limits_{u\rightarrow \infty,\,u\in\mathscr{V}_1^*}F(u)=+\infty,
\end{equation}
which means that $F(u)$ is still proper on $\mathscr{V}_1^*$. Since $F(u(t))'\leq0$,
$u(t)$ lies in a compact region in $\mathscr{V}_1^*$. By Lemma \ref{compact imply CCCPM exist}, the solution exists for $t\in[0,+\infty)$ and there exists a sequence of metrics $r(t_k)$ which converge to $r^*$ as $t_k\uparrow +\infty$.
Hence $F(u(+\infty))=\lim\limits_{k\rightarrow +\infty}F(u(t_k))=F(u^*)=0$.
While Lemma \ref{property of Ricci potential under general condition}
says that $u^*$ is the unique zero and minimum point of $F$, thus we get $u(t)$ converges to $u^*$ as $t\rightarrow +\infty$.
The exponentially convergent rate comes from Theorem \ref{local convergence for NCRF under general conditions}.
\qed

\begin{corollary}\label{global convergence of NCRF}
Suppose $(M, \mathcal{T}, \Phi)$ is a weighted triangulated surface with $\chi(M)\leq 0$.
Then the solution to the flow (\ref{equation of normalized combinatorial Ricci flow on surface}) converges if and only if there exists a constant $R$-curvature metric $r^*$. Furthermore, if the solution converges, it converges exponentially fast to the metric of constant curvature.
\end{corollary}

We expect Corollary \ref{global convergence of NCRF} to be still valid for surfaces with $\chi(M)>0$.
However, things are not so satisfactory for general triangulations.
In fact, things may become very different and complicated.
Assume there exists a constant curvature metric $r^*$, we may compare the Ricci potential (\ref{Ricci potential}) with
\begin{equation*}
G(u)=\int_{u^*}^u\sum_{i=1}^{N}(K_i-K_{i}^*)du_i,
\end{equation*}
where $K_i^*=\frac{2\pi\chi(M)}{\|r^*\|^2}r_i^{*2}$, then
\begin{equation*}
  \begin{aligned}
    F(u)-G(u)
  &=\int_{u_*}^u\sum_{i=1}^{N}(K_i^*-R_{av}r_i^2)du_i\\
  &=2\pi\chi(M)\left(\cfrac{\sum_{i=1}^Nr_i^{*2}(u_i-u_i^*)}{\|r^*\|^2}-\int_{u_*}^u\cfrac{e^{u_1}du_1+\cdots+e^{u_N}du_N}{e^{u_1}+\cdots+e^{u_N}}\right)\\
  &=2\pi\chi(M)\cfrac{\sum_{i=1}^Nr_i^{*2}(u_i-u_i^*)}{\|r^*\|^2}-2\pi\chi(M)\ln(e^{u_1}+\cdots+e^{u_N})\bigg|_{u^*}^u.
  \end{aligned}
\end{equation*}
Restricted to the hypersurface $\mathscr{V}^*\triangleq\{u\in \mathds{R}^N|\sum_i e^{u_i}=\sum_ie^{u_i^*}\}$, the last term becomes zero, hence
\begin{equation*}
F(u)=G(u)+2\pi\chi(M)\cfrac{\sum_{i=1}^Nr_i^{*2}(u_i-u_i^*)}{\|r^*\|^2},\,\,\forall u\in\mathscr{V}^*.
\end{equation*}
By similar arguments in the proof of Theorem \ref{global convergence of NCRF under general condition}, we further have
\begin{equation*}
\lim\limits_{u\rightarrow \infty,\,u\in\mathscr{V}^*}G(u)=+\infty.
\end{equation*}
One may expect that the growth behavior of $F$ is similar to that of $G$.
Unfortunately, there is an obstruction which makes things very complicated. It's easy to see
\begin{equation*}
\lim\limits_{u\rightarrow \infty,\,u\in\mathscr{V}^*}\cfrac{\sum_{i=1}^Nr_i^{*2}(u_i-u_i^*)}{\|r^*\|^2}=-\infty.
\end{equation*}
The growth behavior of $F$ is uncertain unless one can compare the growth rates of these two terms concretely. It seems that we can not expect that $G(u)$ succeeds the last term due to the following example.

\begin{example} \label{example triangluate sphere}
Given a topological sphere, triangulate it into four faces of a single tetrahedron. Fix weights $\Phi\equiv0$. Denote $r^*=\textbf{1}$, then $r^*$ is a constant curvature metric with $R_i\equiv \pi$. By direct calculation, we get $Hess_uF|_{u^*}=(\frac{\sqrt{3}}{6}-\frac{\pi}{4})(4I_4-\textbf{1}\textbf{1}^T)$, which is negative semi-definite with three negative eigenvalues.
Up to scaling, $Hess_u F$ is negative definite at $r^*$. Thus the fixed point $r^*=\textbf{1}$ of the ODE system
$r_i'(t)=\frac{1}{2}(R_{av}-R_i)r_i$ (i.e. the flow (\ref{equation of normalized combinatorial Ricci flow on surface})) is a source. For this particular weighted triangulation, the flow (\ref{equation of normalized combinatorial Ricci flow on surface}) can never converge to the constant curvature metric $r^*$ for any initial value $r(0)$, as $t\rightarrow +\infty$. However, when $t\rightarrow -\infty$, the solution of (\ref{equation of normalized combinatorial Ricci flow on surface}) converges to $r^*$ if $r(0)\in \mathbb{S}^{N-1}(2)$ is close enough to $r^*$.
\end{example}

From the example we know that it is impossible to get an analogous version of Corollary \ref{global convergence of NCRF} when $\chi(M)>0$. However, we can get the following long time existence of (\ref{equation of normalized combinatorial Ricci flow on surface}) by maximum principle.

\begin{theorem} \label{positive Euler number Thm}
Given a weighted triangulated surface $(M, \mathcal{T}, \Phi)$ with $\chi(M)>0$ and
\begin{equation}\label{sphere com-topo condition}
-\sum_{(e,v)\in Lk(I)}(\pi-\Phi(e))+2\pi\chi(F_I)<0, \ \ \forall I:\,\emptyset \subsetneqq I \subsetneqq V.
\end{equation}
Suppose $r_0$ is a circle packing metric with $R_i\geq 0$ for all $i\in V$. Then the solution of (\ref{equation of normalized combinatorial Ricci flow on surface}) with initial metric $r_0$ exists for all time and lies in a compact region in $\mathds{R}^N_{>0}$. Furthermore, there exists $t_n\uparrow +\infty$, such that $r(t_n)$ converges to a constant curvature metric $r^*$.
\end{theorem}

\textbf{Proof.}
By Lemma \ref{compact imply CCCPM exist}, we just need to prove that $\{r(t)\}\subset\subset \mathds{R}^N_{>0}$.
We prove it by contradiction. Suppose
$$\overline{\{r(t)\}}\cap \partial (\mathds{R}^N_{>0}\cap \mathbb{S}^{N-1})\neq \emptyset,$$
then there exists a sequence $t_n$ such that
$$r(t_n)\rightarrow r^*\in \partial (\mathds{R}^N_{>0}\cap \mathbb{S}^{N-1}),$$
where $r^*=(0,\cdots, 0, r^*_{|I|+1}, \cdots, r_N^*)$ for some nonempty proper subset $I$ of  $V$.
By Proposition 4.1 of \cite{CL1}, we have
\begin{equation}\label{contradiction}
\lim_{n\rightarrow +\infty}\sum_{i\in I}K_i(r(t_n))=-\sum_{(e,v)\in Lk(I)}(\pi-\Phi(e))+2\pi\chi(F_I)<0.
\end{equation}
However, by Corollary \ref{preserve positive curv}, $R_i\geq 0$ is preserved along the flow
(\ref{equation of normalized combinatorial Ricci flow on surface}). Thus $K_i\geq 0$ for all $i\in V$,
which implies $\lim_{n\rightarrow +\infty}\sum_{i\in I}K_i(r(t_n))\geq 0$. This contradicts (\ref{contradiction}).
\qed

\begin{remark}
It is easy to see that the triangulation of sphere in Example \ref{example triangluate sphere} does not satisfy the condition (\ref{sphere com-topo condition}).
\end{remark}

\subsection{The prescribing curvature problem}
Using modified combinatorial Ricci flow, we can consider the prescribing curvature problem
on a weighted triangulated surface $(M, \mathcal{T}, \Phi)$.

\begin{definition}\label{definition of modified CRF}
Suppose $(M, \mathcal{T}, \Phi)$ is a weighted triangulated surface with circle packing metric $r$, $\overline{R}\in C(V)$
is a function defined on $M$. The modified combinatorial Ricci flow with respect to $\overline{R}$ is defined to be
\begin{equation}\label{equation of modified CRF}
\frac{dg_i}{dt}=(\overline{R}_i-R_i)g_i,
\end{equation}
where $g_i=r_i^2$ as before.
\end{definition}
Note that the total measure $\mu(M)=||r||^2$ may vary along the modified combinatorial Ricci flow (\ref{equation of modified CRF}),
which is different from that of (\ref{equation of normalized combinatorial Ricci flow on surface}).

$\overline{R}$ is called admissible if there is a circle packing metric $\overline{r}$ with curvature $\overline{R}$. For a given
function $\overline{R}\in C(V)$, we can introduce the following modified Ricci potential
\begin{equation}\label{definition of modified Ricci potential}
\overline{F}(u)=\int_{u_0}^u\sum_{i=1}^N\left(K_i-\overline{R}_ir_i^2\right)du_i.
\end{equation}
It is easy to check that the modified Ricci potential $\overline{F}$ is well-defined. Furthermore, by direct
calculations, we have
\begin{equation*}
\begin{aligned}
Hess_u\overline{F}=L-\Sigma^{\frac{1}{2}}
\left(
       \begin{array}{ccc}
       \overline{R}_1 &   &   \\
        & \ddots &   \\
        &   & \overline{R}_N \\
       \end{array}
\right)\Sigma^{\frac{1}{2}}.
\end{aligned}
\end{equation*}

The following lemma is useful.
\begin{lemma}(\cite{CL1})\label{lemma of injective}
Suppose $\Omega\subset \mathds{R}^N$ is convex, the function $h:\Omega\rightarrow \mathds{R}$
is strictly convex, then the map $\nabla h:\Omega\rightarrow \mathds{R}^N$ is injective.
\end{lemma}

\begin{theorem}\label{existence for prescribing curvature problem}
Suppose $(M, \mathcal{T}, \Phi)$ is a weighted triangulated surface and $\overline{R}\in C(V)$ is a function defined on $M$.
\begin{itemize}
  \item[(1)] If the solution to the modified flow (\ref{equation of modified CRF}) converges, then $\overline{R}$ is admissible.
  \item[(2)] If $\overline{R}_i\leq 0$ for all $i$, but not identically zero, and $\overline{R}$ is admissible by a metric $\overline{r}$. Then $\overline{r}$ is the unique metric in $\mathds{R}^N_{>0}$ such that its curvature is $\overline{R}$. Moreover, the solution to the modified flow (\ref{equation of modified CRF}) converges exponentially fast to $\overline{r}$.
\end{itemize}
\end{theorem}

\textbf{Proof.}
The first part is obviously, and $\overline{R}$ is admissible by metric $r(+\infty)$. For the second part, notice that
\begin{equation*}
Hess_u\overline{F}=L-\Sigma^{\frac{1}{2}}diag\{\overline{R}_1,\cdots,\overline{R}_N\}\Sigma^{\frac{1}{2}},
\end{equation*}
it is easy to check that, if $\overline{R}_i\leq 0$ for $i=1, \cdots, N$ and not identically zero, $Hess_u\overline{F}$ is
positive definite. By Lemma \ref{lemma of injective}, $\nabla_u \overline{F}=(K_1-\overline{R}_1r_1^2,\cdots,K_N-\overline{R}_Nr_N^2)^T$
is an injective map from $u\in \mathds{R}^N$ to $\mathds{R}^N$. Hence $\overline{r}$ is the unique zero point of $\nabla_u \overline{F}$. This fact implies that $\overline{r}$ is the unique metric in $\mathds{R}^N_{>0}$ such that it's curvature is $\overline{R}$. By Lemma B.1 in \cite{Ge}, we know that $\overline{F}$ is proper and $\lim\limits_{u\rightarrow\infty}\overline{F}(u)=+\infty$. Furthermore, $\frac{d}{dt}F(u(t))=-\sum_ir_i^{-2}(K_i-\overline{R}_ir_i^2)^2\leq 0$ implies that the solution to (\ref{equation of modified CRF}) lies in a compact region. The rest of the proof is the same as that of Theorem \ref{global convergence of NCRF under general condition}, so we omit it here. \qed
\begin{remark}
Theorem \ref{existence for prescribing curvature problem} could be taken as a discrete version of a result obtained by Kazdan and Warner in \cite{KW1}.
\end{remark}

\begin{remark}\label{prescribing problem for R=0}
The second part of Theorem \ref{existence for prescribing curvature problem} implies that $\chi(M)<0$. If $\overline{R}_i=0$ for all $i$, the corresponding prescribing curvature problem is already solved in Corollary \ref{global convergence of NCRF}. In this case, the metric $\overline{r}$ is not unique. However, it is unique up to scaling. This is slightly different from Theorem \ref{existence for prescribing curvature problem}.
\end{remark}

\subsection{Existence and uniqueness of constant curvature metric}
\label{section-existence and uniqueness of cccpm}
We first condider the uniqueness of constant $R$-curvature metric.
\begin{theorem}\label{uniqueness of CCM}
Suppose $(M, \mathcal{T}, \Phi)$ is a weighted triangulated surface. $c^*\in\mathds{R}$ is a constant.
\begin{itemize}
  \item[(1)] If $c^*<0$, there exists at most one circle packing metric with curvature $R=c^*\textbf{1}$.
  \item[(2)] If $c^*=0$, there exists at most one metric with curvature $R=c^*\textbf{1}$ up to scaling.
\end{itemize}
\end{theorem}
\textbf{Proof.}
The first part is just a corollary of the second part of Theorem \ref{existence for prescribing curvature problem}. The second part is already proved by Thurston \cite{T1}.
\qed \\

If $\chi(M)\leq0$, Theorem \ref{uniqueness of CCM} implies that the constant curvature metric (if exists) is unique or unique up to scaling.
If $\chi(M)>0$, there may be several different constant curvature metrics on $\mathbb{S}^{N-1}$ and we have the following example.
\begin{example}\label{example-not unique cccpm}
Consider the weighted triangulation of the sphere in Example \ref{example triangluate sphere} with tangential circle packing metric.
Let $r_1=1$, $r_i=x>0$, $2\leq i\leq4$ and then $x=\frac{\cos\theta}{1-\cos\theta}$,
where $\theta$ is the inner angle of $\triangle v_1v_2v_3$ at the vertex $v_3$. \\
\begin{figure}[hbp]
\center
\begin{minipage}[t]{0.75\linewidth}
\centering
\includegraphics[width=0.55\textwidth]{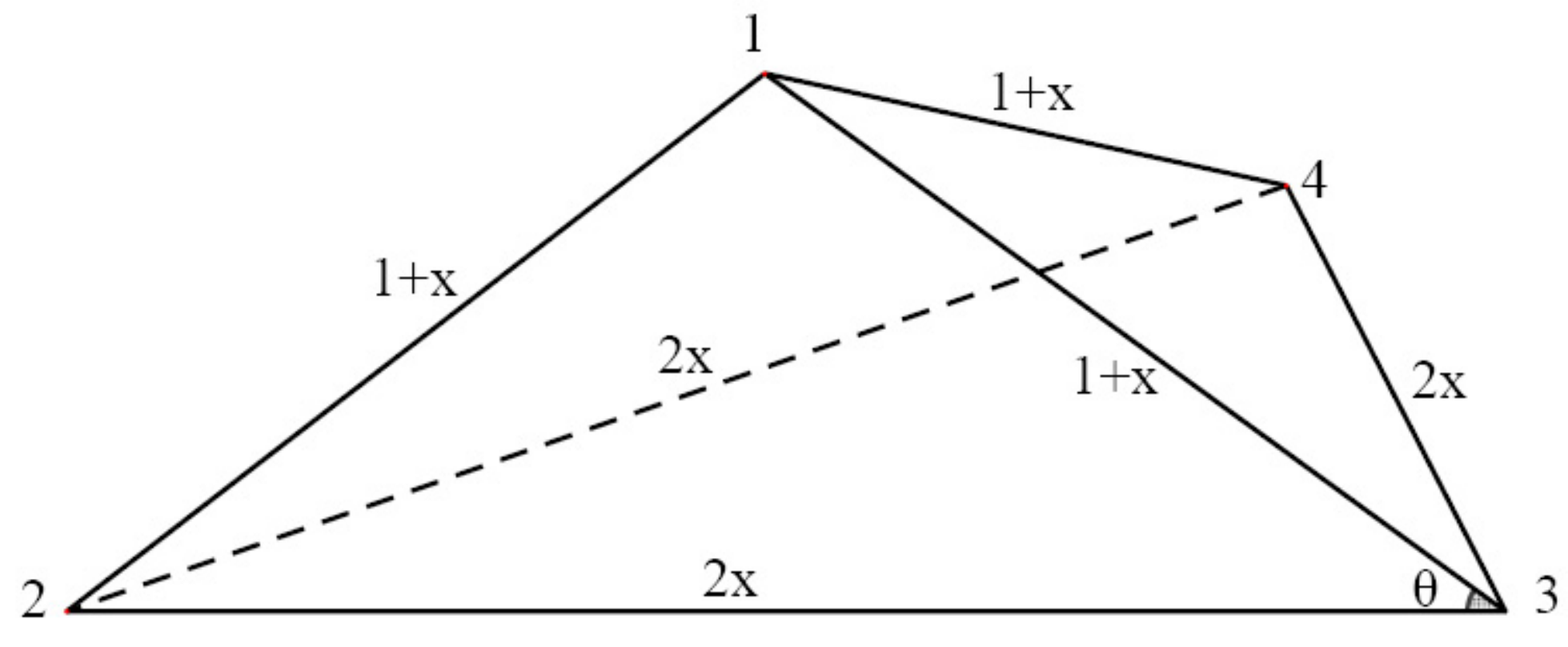}
\caption{{tetrahedron triangulation of sphere}}
\label{tetrahedron triangulation of sphere}
\end{minipage}
\end{figure}
\\
By direct calculations, $K_1=6\theta-\pi$, $K_i=\frac{5}{3}\pi-2\theta$, $2\leq i\leq4$. Constant $R$-curvature implies that
$$\frac{\frac{5}{3}\pi-2\theta}{x^2}=\frac{6\theta-\pi}{1^2}.$$
Using MATLAB, we can solve the equation and get two solutions $x=1$ and $x\approx 5.9487$. Thus we find at least two constant curvature metrics $r_1^*=\textbf{1}$ and $r_2^*\approx(1, 5.9487, 5.9487, 5.9487)$.
\end{example}

In the following of this subsection, we consider the existence of constant $R$-curvature metric.
For the classical discrete Gauss curvature $K_i$,
Theorem \ref{Thm-thurston} states that the existence of constant $K$-curvature metric is equivalent to
\begin{equation}\label{Thuston's condition}
2\pi\chi(M)\frac{|I|}{|V|}>-\sum_{(e,v)\in Lk(I)}(\pi-\Phi(e))+2\pi\chi(F_I)
\end{equation}
for any nonempty proper subset $I$ of $V$.
Theorem \ref{Thm-thurston-MardenRodin-ChowLuo} further states that the space of all admissible classical discrete Gauss curvature $K$ is
\begin{equation*}
\mathscr{Y}\triangleq\mathscr{K}_{GB} \cap \Big(\mathop{\cap} \limits_{\phi\,\neq\,I \,\subsetneqq\,V} \mathscr{Y}_I \Big),
\end{equation*}
where $\mathscr{K}_{GB}$ and $\mathscr{Y}_I$ are defined by (\ref{Gauss-Bonnet plane}) and (\ref{Y_I half space}) respectively.
If $r^*$ determines a constant $R$-curvature, then the $K$-curvature $K^*=(K_1^*,\cdots,K_N^*)^T$ is admissible, where $K_i^*=\frac{2\pi\chi(M)}{||r^*||^{2}_{2}}r_i^{*2}$. Hence $K^*\in\mathscr{Y}_I$, i.e.
\begin{equation}\label{the comb topo condition of CL}
2\pi\chi(M)\frac{\sum_{i\in I}r_i^{*2}}{||r^*||^{2}_{2}}>-\sum_{(e,v)\in Lk(I)}(\pi-\Phi(e))+2\pi\chi(F_I)
\end{equation}
for any nonempty proper subset $I$ of $V$.
(\ref{the comb topo condition of CL}) is necessary for the existence of constant $R$-curvature metric,
which involves the combinatorial, topological and metric structure of the triangulated surface.
In the case of $\chi(M)\leq 0$, the combinatorial-topological-metric condition (\ref{the comb topo condition of CL}) is also sufficient.
In fact, we have the following result.

\begin{theorem}\label{Thm-combtopo-condition-CCCPmetric}
Suppose $(M, \mathcal{T}, \Phi)$ is a weighted triangulated surface with $\chi(M)\leq 0$.
Then there exists a constant $R$-curvature metric if and only if
there exists a circle packing metric $r^*$ such that, for any nonempty proper subset $I$ of vertices $V$,
(\ref{the comb topo condition of CL}) is valid.
\end{theorem}
\textbf{Proof.}
In the case of $\chi(M)=0$, a zero $R$-curvature metric is exactly a zero $K$-curvature metric.
The conclusion is contained in Theorem \ref{Thm-thurston}.

In the case of $\chi(M)<0$, We just need to prove the sufficient part.
Set $K_i^*=2\pi\chi(M)\frac{\sum_{i\in I}r_i^{*2}}{\|r^*\|^{2}_{2}}$.
By Theorem \ref{Thm-thurston-MardenRodin-ChowLuo},
(\ref{the comb topo condition of CL}) implies that $K^*=(K_1^*,\cdots,K_N^*)^T$ is an admissible $K$-curvature.
Suppose $K(r_0)=K^*$. Then $r_0$ is a circle packing metric with $R_i(r_0)<0$ for all $i\in V$.
By Theorem \ref{existence theorem for negative curvature},
the combinatorial Ricci flow (\ref{equation of normalized combinatorial Ricci flow on surface}) with $r_0$ as the
initial metric converges to a negative constant $R$-curvature metric.
Therefore, there exists a constant $R$-curvature metric on $(M, \mathcal{T}, \Phi)$. \qed

The key point of the proof of Theorem \ref{Thm-combtopo-condition-CCCPmetric} is to find a circle packing metric
with all $R$-curvature negative and then use the combinatorial Ricci flow (\ref{equation of normalized combinatorial Ricci flow on surface})
to evolve it to a constant $R$-curvature metric. Therefore, we can generalize Theorem \ref{Thm-combtopo-condition-CCCPmetric}
to the following form.

\begin{theorem}\label{main existence thm}
Suppose $(M, \mathcal{T}, \Phi)$ is a weighted triangulated surface.
\begin{itemize}
  \item[(1)] There exists a negative constant curvature metric if and only if $\mathscr{Y}\cap\mathds{R}^N_{<0}\neq\phi$.
  \item[(2)] There exists a zero curvature metric if and only if $0\in\mathscr{Y}$.
  \item[(3)] If there exists a positive constant curvature metric, then $\mathscr{Y}\cap\mathds{R}^N_{>0}\neq\phi$.
\end{itemize}
\end{theorem}
We include the case of $\chi(M)>0$ in Theorem \ref{main existence thm}. In this case, if there exists a constant $R$-curvature metric, then
there exists a circle packing metric $r^*$ satisfying (\ref{the comb topo condition of CL}), which
implies $\mathscr{Y}\cap\mathds{R}^N_{>0}\neq\phi$.

Similar to Thurston's condition (\ref{condition-Thurston}),
conditions in Theorem (\ref{main existence thm}) also show that the combinatorial structure of the triangulation and the topology of surface,
which have no relation to the circle packing metrics, contain some information of $R$-curvature.

For the surfaces with $\chi(M)>0$, the property of constant $R$-curvature metric seems complicated and confusing as exhibited in Example \ref{example triangluate sphere} and Example \ref{example-not unique cccpm}.
We do not know whether the condition $\mathscr{Y}\cap\mathds{R}^N_{>0}\neq\phi$ is sufficient for
the existence of constant $R$-curvature metric.
Even through, we can get some sufficient conditions for the existence of constant $R$-curvature metric by the discrete maximum principle.

\begin{theorem}\label{general existence theorem}
Suppose $(M, \mathcal{T}, \Phi)$ is a weighted triangulated surface. If there is a circle packing metric $r$ satisfying $R_{i}\geq0$ for all $i\in V$, and
\begin{equation*}
-\sum_{(e,v)\in Lk(I)}(\pi-\Phi(e))+2\pi\chi(F_I)<0
\end{equation*}
for any nonempty subset $I$ of $V$, then there exists a non-negative constant $R$-curvature metric on $(M, \mathcal{T}, \Phi)$.
\end{theorem}

\textbf{Proof.}
By the maximum principle, i.e. Theorem \ref{Maximum priciple},
$R_{i}\geq 0$ is preserved along the normalized combinatorial Ricci flow (\ref{equation of normalized combinatorial Ricci flow on surface}).
Using similar arguments in the proof of Theorem \ref{positive Euler number Thm}, we get the proof.
\qed

\section{Combinatorial Calabi flow on surfaces}\label{Calabi flow section}

Given a compact complex manifold admitting at least one K\"{a}hler metric, to find the extreme metric
which minimizes the $L^2$ norm of the curvature tensor in a given principal cohomology class,
Calabi \cite{Ca} introduced the Calabi flow, which could be written as
\begin{equation}\label{Calabi flow}
\frac{\partial g}{\partial t}=\Delta_g K\cdot g
\end{equation}
on a Riemannian surface. Chru\'{s}ciel \cite{CP} proved that the Calabi flow (\ref{Calabi flow}) exists
for all time and
converges to a constant Gaussian curvature metric on closed surfaces using the Bondi mass estimate,
assuming the existence of the constant
Gaussian curvature in the background which is ensured by the the uniformization theorem. Chang \cite{CSC} pointed
out that Chru\'{s}ciel's results still hold for arbitrary initial metric, which implies the
uniformization theorem on closed surfaces with genus greater than one.
Chen \cite{CXX} gave a geometrical proof of the long-time existence and the convergence of the Calabi flow on closed surfaces.
The first author \cite{Ge} first introduced the notion of combinatorial Calabi flow in Euclidean background geometry
and proved that the convergence
of the flow is equivalent to the existence of constant classical discrete Gauss curvature.
Then we \cite{GX1} studied the combinatorial Calabi flow in hyperbolic background geometry.
We also use the combinatorial Calabi flow to study some constant combinatorial curvature problem on
3-dimensional triangulated manifolds \cite{GX2}.
In this section, we study the $R$-curvature problem by combinatorial Calabi flow.

\begin{definition}\label{combinatorial Calabi flow}
For a weighted triangulated surface $(M, \mathcal{T}, \Phi)$ with circle packing metric $r$,
the combinatorial Calabi flow is defined as
\begin{equation}\label{equation of CCF}
\frac{dg_i}{dt}=\Delta R_i\cdot g_i,
\end{equation}
where $g_i=r_i^2$ and $\Delta$ is the Laplacian operator given by (\ref{definition of Laplacian}).
\end{definition}

It is easy to check that the total measure $\mu(M)=\|r\|^2$ of $M$ is invariant along the combinatorial Calabi flow (\ref{equation of CCF}).
Interestingly, the combinatorial curvature $R_i$ evolves according to
$$\frac{dR_i}{dt}=-R_i\Delta R_i-\Delta^2R_i,$$
which has almost the same form as that of the scalar curvature along the smooth Calabi flow on surfaces.
We can rewrite the combinatorial Calabi flow (\ref{equation of CCF}) as
$$\frac{dr_i}{dt}=-\frac{1}{2r_i}\sum_j\frac{\partial K_i}{\partial u_j}R_j=-\frac{1}{2r_i}\sum_j\frac{\partial K_i}{\partial u_j}(R_j-R_{av}).$$
Set $$\Gamma(r)_i=-\frac{1}{2r_i}\sum_j\frac{\partial K_i}{\partial u_j}(R_j-R_{av}).$$
If there exists a constant curvature metric $r^*$, then we have
\begin{equation*}
\begin{aligned}
D_{j}\Gamma_i|_{r^*}
=-\frac{1}{r_ir_j}\left(\sum_k\frac{1}{r_k^2}\frac{\partial K_i}{\partial u_k}\frac{\partial K_j}{\partial u_k}-R_{av}\frac{\partial K_i}{\partial u_j}\right).
\end{aligned}
\end{equation*}
If we further assume that the condition (\ref{lambda(lambda)>c general condition}), i.e. $\lambda_1(-\Delta)>R_{av}^*=\frac{2\pi\chi(M)}{||r^*||^2}$, is satisfied, then $D\Gamma|_{r^*}$ is negative semi-definite with rank $N-1$ and kernel $\{t\textbf{1}|t\in \mathds{R}\}$. From this fact we can get local convergence results similar to Theorem \ref{local convergence for NCRF under general conditions}, Corollary \ref{local convergence for NCRF}. We can also get results similar to Lemma \ref{compact imply CCCPM exist}, Corollary \ref{convergence imply CCCPM exist}, Theorem \ref{global convergence of NCRF under general condition} and Corollary \ref{global convergence of NCRF}. We just state the following main theorem for this section here.

\begin{theorem}\label{global convergence of CCF}
Suppose $(M, \mathcal{T}, \Phi)$ is a weighted triangulated surface with $\chi(M)\leq 0$,
then the combinatorial Calabi flow (\ref{equation of CCF}) converges if and only if there exists a constant $R$-curvature metric $r^*$.
\end{theorem}

\textbf{Proof.}
If the solution of (\ref{equation of CCF}) converges, $u(+\infty)$ must be a critical point of this ODE system, i.e. $\Delta R(+\infty)=-\Sigma^{-1}LR(+\infty)=0$, which implies that $R(+\infty)$ belongs to the kernel of $L$
and hence $r(+\infty)$ is a constant curvature metric.

Assume there exists a constant curvature metric $r^*$. Write the flow (\ref{equation of CCF}) in the matrix form $\dot{u}=\Delta R=-\Sigma^{-1}LR$ or $\dot{r^2}=-LR$, where $r^2=(r_1^2,\cdots,r_N^2)^T$. It is easy to see that
$$\frac{d}{dt}F(u(t))=-(K-R_{av}r^2)^T\Sigma^{-1}LR=-R^TLR\leq 0.$$
Using the properties of the Ricci potential, we can get the convergence result by similar arguments
in the proof of Theorem \ref{global convergence of NCRF under general condition}.\qed\\

Analogous to the combinatorial Ricci flow, we can also use the combinatorial Calabi flow to study the combinatorial
prescribing curvature problem. In fact, we have the following result.

\begin{theorem}\label{existence for prescribing curvature problem using CCF}
Suppose $(M, \mathcal{T}, \Phi)$ is a weighted triangulated surface, $\overline{R}\in C(V)$ is a function defined on $M$ with $\overline{R}_i\leq 0$ for all $i$. Then $\overline{R}$ is admissible if and only if the solution of the modified combinatorial Calabi flow
\begin{equation}
\frac{dg_i}{dt}=\Delta (R-\overline{R})_ig_i
\end{equation}
exists for all time and converges.\qed
\end{theorem}

The proof is just a combination of that of Theorem \ref{existence for prescribing curvature problem}, Remark \ref{prescribing problem for R=0}
and Theorem \ref{global convergence of CCF}, we omit it here.

Following \cite{Ge}, we can also introduce a notion of energy to study the combinatorial curvature $R_i$.

\begin{definition}
For a weighted triangulated surface $(M, \mathcal{T}, \Phi)$ with circle packing metric $r$,
the combinatorial Calabi energy is defined as
\begin{equation}
\widetilde{C}(r)=\sum_{i=1}^{N}\varphi_i^2,
\end{equation}
where $\varphi_i=K_i-\frac{2\pi \chi_M}{||r||^2}r_i^2$.
\end{definition}

Consider the combinatorial Calabi energy $\widetilde{C}$ as a function of $u$, we have
$\nabla_u\widetilde{C}=2A^T\varphi$, where
$$A=\frac{\partial (\varphi_1, \cdots, \varphi_N)}{\partial (u_1, \cdots, u_N)}
=\left(
   \begin{array}{ccc}
     \frac{\partial \varphi_1}{\partial u_1} & \cdots & \frac{\partial \varphi_1}{\partial u_N} \\
     \vdots & \ddots & \vdots \\
     \frac{\partial \varphi_N}{\partial u_1} & \cdots & \frac{\partial \varphi_N}{\partial u_N} \\
   \end{array}
 \right).
$$
Note that $A$ is in fact given by (\ref{Hession of F}). Thus, if $\chi(M)\leq 0$, $A$ is symmetric,  positive semi-definite with
rank $N-1$ and kernel $\{t \mathbf{1}|t\in \mathds{R}\}$.
Using this fact, we can define a new flow, which is the gradient flow of $\widetilde{C}$. We call it the modified
combinatorial Calabi flow.

\begin{definition}
For a weighted triangulated surface $(M, \mathcal{T}, \Phi)$ with circle packing metric $r$,
the modified combinatorial Calabi flow is defined as
\begin{equation}\label{2-dim CFF with energy}
\dot{u}=-\frac{1}{2}\nabla_u\widetilde{C},
\end{equation}
or equivalently,
\begin{equation}
\dot{u}=-A^T\varphi,
\end{equation}
where $\varphi=(\varphi_1, \cdots, \varphi_N)^T$.
\end{definition}

Note that $\sum_i\varphi_i=0$, so in the case of $\chi(M)\leq 0$, if the flow (\ref{2-dim CFF with energy}) converges, it converges
to the constant curvature metric.
Along the flow (\ref{2-dim CFF with energy}), we have
$$\dot{\varphi}=-AA^T\varphi,\ \dot{\widetilde{C}}=-2||A^T\varphi||^2\leq 0$$
and
$$\frac{d}{dt}F(u(t))=(\nabla_uF)^T\cdot \dot{u}=\varphi^T\cdot \dot{u}=-\varphi^TA^T\varphi\leq 0,$$
which implies that $F(u(t))$ is decreasing along the flow (\ref{2-dim CFF with energy}).

Following the arguments in the proof of Theorem \ref{global convergence of NCRF under general condition} and \ref{global convergence of CCF},
we can derive the following result.

\begin{theorem}\label{global convergence of modified CCF}
Suppose $(M, \mathcal{T}, \Phi)$ is a weighted triangulated surface with $\chi(M)\leq 0$, then the existence of constant curvature metric $r^*$ is equivalent to the convergence of the modified combinatorial Calabi flow (\ref{2-dim CFF with energy}).\qed
\end{theorem}

\section{Combinatorial $\alpha$-curvature and combinatorial $\alpha$-flows} \label{section-alpha}
In the previous sections, we investigated the properties of the new discrete Gauss curvature $R_i=K_i/r_i^2$.
Since the area of the disk packed at $i$ is just $\pi r_i^2$, the denominator of the curvature $r_i^2$ (omitting the efficient $\pi$) may also be considered as an ``area element" attached to vertex $i$. We want to know whether there exists other types of ``area element". Suppose $A_i$ is a general ``area element", which is an analogy of
the volume element in the smooth case, the combinatorial Gauss curvature could be defined as $R_i=K_i/A_i$. The average curvature should be $2\pi\chi(M)/\sum A_i$, which is an analogy of $\int R d\mu/\int d\mu$ in the smooth case. If we expect the functional
$$F(u)=\int_{u_0}^u\sum_{i=1}^{N}(K_i-\frac{2\pi\chi(M)}{\sum A_i}A_i)du_i$$
to be well defined, the following formula
\begin{equation}\label{area elment equation}
\frac{\partial}{\partial u_i}\left(\frac{A_j}{\sum_k A_k}\right)=\frac{\partial}{\partial u_j}\left(\frac{A_i}{\sum_k A_k}\right)
\end{equation}
should be satisfied for all $i$ and $j$. It is a trivial observation that $A_i=r_i^{\alpha}$ always satisfies (\ref{area elment equation}) for all $\alpha\in \mathds{R}$.

When $\chi(M)\leq0$, we have already seen in the previous sections that the combinatorial Ricci flow and Calabi flow are good enough to evolve circle packing metrics to a constant $R$-curvature metric. However, in the case of $\chi(M)>0$, we do not know how to evolve it. Maybe this is because that the growth behavior of $F$ is uncertain (in integral level), or equivalently, $Hess_u F$ is not positive semi-definite (in differential level), even for the simplest triangulation of $\mathbb{S}^2$ (see Example \ref{example triangluate sphere}).

The above two reasons motivate us to consider a new ``area element" $A_i=r_i^{\alpha}$. We find that, for $\chi(M)>0$, if $\alpha<0$, combinatorial flow methods are good enough to evolve the curvature to a constant. For $\chi(M)\leq0$, if $\alpha>0$, the properties of combinatorial flows are almost the same with previous sections. It is very interesting that \cite{CL1} and \cite{Ge} can be included into the case of $\alpha=0$.

\begin{definition}\label{definition of alpha curvature}
For a weighted triangulated surface $(M, \mathcal{T}, \Phi)$ with a circle packing metric $r$,
the combinatorial $\alpha$-curvature is defined as
\begin{equation}\label{alpha curvature}
R_{\alpha,i}=\frac{K_i}{r_i^\alpha},
\end{equation}
where $\alpha$ is a real number.
\end{definition}

For this type of curvature, we can also consider the corresponding constant curvature problem and prescribing curvature problem.
As the methods are all the same as that we dealt with $R$-curvature, we will give only the outline and skip the details in the following.

The measure defined on $M$ is now $\mu_\alpha(i)=r_i^\alpha$ and the average curvature is now
$$R_{\alpha, av}=\frac{\int_M R_{\alpha}d\mu_\alpha}{\int_M d\mu_\alpha}=\frac{2\pi\chi(M)}{||r||_\alpha^\alpha},$$
where $||r||_\alpha=(\sum_{i=1}^Nr_i^\alpha)^{\frac{1}{\alpha}}$ and $||r||_\alpha^\alpha$ is defined to be $N$ in the case of $\alpha=0$.
In this section, we set $u=\ln r_i$ and define the Ricci potential as
$$F_\alpha(u)=\int_{u_0}^u\sum_{i=1}^N(K_i-R_{\alpha,av}r_i^\alpha)du_i,$$
where $u_0\in \mathds{R}^N$ is an arbitrary point. By direct calculations, it is easy to check that the Ricci potential
$F_\alpha$ is well-defined. We further have
$\nabla_uF_\alpha=K-R_{\alpha,av}r^\alpha$ and
\begin{equation*}
\begin{aligned}
Hess_uF_\alpha
&=\widetilde{L}-\alpha R_{\alpha,av}\Sigma^{\frac{\alpha}{4}}\left(I-\frac{r^{\frac{\alpha}{2}}\cdot (r^{\frac{\alpha}{2}})^T}{||r||_\alpha^\alpha}\right)\Sigma^{\frac{\alpha}{4}}\\
&=\left(
   \begin{array}{ccc}
     r_1^{\frac{\alpha}{2}} &   &   \\
       & \ddots &   \\
       &   & r_1^{\frac{\alpha}{2}} \\
   \end{array}
 \right)
 \left(\Lambda_\alpha-\alpha R_{\alpha,av}\left(I-\frac{r^{\frac{\alpha}{2}}\cdot (r^{\frac{\alpha}{2}})^T}{||r||_\alpha^\alpha}\right)\right)
\left(
   \begin{array}{ccc}
     r_1^{\frac{\alpha}{2}} &   &   \\
       & \ddots &   \\
       &   & r_1^{\frac{\alpha}{2}} \\
   \end{array}
 \right),
\end{aligned}
\end{equation*}
where $r^\alpha=(r_1^\alpha, \cdots, r_N^\alpha)^T$,
$\Lambda_\alpha=\Sigma^{-\frac{\alpha}{4}}\widetilde{L}\Sigma^{-\frac{\alpha}{4}}$
and $\widetilde{L}_{ij}=\frac{\partial K_i}{\partial u_j}$
with $u_j=\ln r_j$. Note that the matrix
$I-\frac{r^{\frac{\alpha}{2}}\cdot (r^{\frac{\alpha}{2}})^T}{||r||_\alpha^\alpha}$
has eigenvalues 1 ($N-1$ times) and 0 (1 time) and kernel $\{cr^{\frac{\alpha}{2}}| c\in \mathds{R}\}$.
Following the arguments in the proof of Lemma \ref{property of Ricci potential under general condition},
we have, if the first positive eigenvalue of $\Lambda_\alpha$ satisfies
$$\lambda_1(\Lambda_\alpha)> \alpha R_{\alpha,av},$$
$Hess_uF_\alpha$ is positive semi-definite with rank $N-1$ and kernel $\{c\textbf{1}| c\in \mathds{R}\}$.
Especially, if $\alpha\chi(M)\leq 0$, $Hess_uF_\alpha$ is positive semi-definite with rank $N-1$
and kernel $\{c\textbf{1}| c\in \mathds{R}\}$.

We can also define the following modification of the combinatorial Ricci flow, which is called $\alpha$-Ricci flow.

\begin{definition}
For a weighted triangulated surface $(M, \mathcal{T}, \Phi)$ with circle packing metric $r$, the
$\alpha$-Ricci flow is defined to be
\begin{equation}\label{alpha flow}
\begin{aligned}
\frac{dr_i}{dt}=-R_{\alpha,i}r_i.
\end{aligned}
\end{equation}
\end{definition}

The normalization of the $\alpha$-Ricci flow is
\begin{equation}\label{normalized alpha flow}
\begin{aligned}
\frac{dr_i}{dt}=(R_{\alpha,av}-R_{\alpha,i})r_i.
\end{aligned}
\end{equation}
Notice that, when $\alpha=0$, the flow (\ref{alpha flow}) and the normalized flow (\ref{normalized alpha flow})
are just Chow and Luo's combinatorial Ricci flows \cite{CL1}.
In this case, $\Pi_{i=1}^N r_i$ (or $\sum_{i=1}^N u_i$) is invariant along the normalized flow (\ref{normalized alpha flow}). When $\alpha\neq0$, $\|r\|_\alpha$ (or $\sum_{i=1}^N e^{\alpha u_i}$) is invariant along the normalized flow (\ref{normalized alpha flow}).
Along the normalized $\alpha$-Ricci flow, the $\alpha$-curvature evolves according to
\begin{equation}\label{evolution of curv under alpha Ricci flow}
\begin{aligned}
\frac{dR_{\alpha, i}}{dt}=\frac{1}{r_i^\alpha}\sum_{j=1}^N\left(-\frac{\partial K_i}{\partial u_j}\right)R_{\alpha, j}+\alpha R_{\alpha,i}(R_{\alpha,i}-R_{\alpha,av}).
\end{aligned}
\end{equation}
For $\alpha=2$, this is just the evolution equation (\ref{evolution of R_i under NCRF}) derived in Lemma \ref{lemma for evolution of R_i under NCRF}.
As the evolution equation (\ref{evolution of curv under alpha Ricci flow}) is still a heat-type equation, the discrete maximum principle, i.e. Theorem \ref{Maximum priciple}, could be applied to this equation.
The evolution of curvature (\ref{evolution of curv under alpha Ricci flow}) under the normalized
$\alpha$-Ricci flow suggests us to define the $\alpha$-Laplacian as
\begin{equation}\label{definition of alpha Laplacian}
\Delta_\alpha f_i
=\frac{1}{r_i^\alpha}\sum_{j=1}^N \left(-\frac{\partial K_i}{\partial u_j}\right) f_j
=\frac{1}{r_i^\alpha}\sum_{j\sim i} \left(-\frac{\partial K_i}{\partial u_j}\right) (f_j-f_i),
\end{equation}
where $f\in C(V)$ and $u_i=\ln r_i$. Using this Laplacian, we can define the combinatorial $\alpha$-Calabi flow as
\begin{equation}\label{definition of alpha Calabi flow}
\frac{dr_i}{dt}=(\Delta_\alpha R_{\alpha})_ir_i.
\end{equation}
Similarly, denote $\varphi_{\alpha,i}=K_i-\frac{2\pi \chi(M)}{||r||^{\alpha}_{\alpha}}r_i^{\alpha}$ and $\widetilde{C}_{\alpha}(r)=\sum_{i=1}^{N}\varphi_{\alpha,i}^2$, we can define a modified $\alpha$-Calabi flow as
\begin{equation}\label{modified alpha Calabi flow}
\dot{u}=-\frac{1}{2}\nabla_u\widetilde{C}_{\alpha}.
\end{equation}

Applying the discrete maximal principle to $\alpha$-Ricci flow (\ref{normalized alpha flow}),
we can get similar existence results for constant $\alpha$-curvature metric.
\begin{theorem}\label{alpha-converge thm for negative initial curvature}
Suppose $(M, \mathcal{T}, \Phi)$ is a weighted triangulated surface. If the initial metric $r(0)$ satisfying $\alpha R_{\alpha,i}(0)<0$ for all $i\in V$, then the normalized $\alpha$-Ricci flow (\ref{normalized alpha flow}) converges to a constant $\alpha$-curvature metric.
\end{theorem}

\textbf{Proof.}
Note that the $\alpha$-curvature $R_{i,\alpha}$ evolves according to (\ref{evolution of curv under alpha Ricci flow})
along the normalized $\alpha$-Ricci flow (\ref{normalized alpha flow}). The maximum principle, i.e. Theorem \ref{Maximum priciple},
is valid for this equation. By the maximum principle, if $\alpha>0$ and $R_{\alpha, i}<0$ for all $i\in V$, we have
$$\left(R_{\alpha, \min}(0)-R_{\alpha, av}\right)e^{\alpha R_{\alpha, av}t}\leq R_{\alpha,i}-R_{\alpha, av}\leq
R_{\alpha, av}(1-\frac{R_{\alpha, av}}{R_{\alpha, \max}(0)})e^{\alpha R_{\alpha, av}t}.$$
If $\alpha<0$ and $R_{\alpha, i}>0$ for all $i\in V$, we have
$$\frac{R_{\alpha,av}}{R_{\alpha,\min}(0)}\left(R_{\alpha, \min}(0)-R_{\alpha, av}\right)e^{\alpha R_{\alpha, av}t}\leq R_{\alpha,i}-R_{\alpha, av}\leq
(R_{\alpha, \max}(0)-R_{\alpha, av})e^{\alpha R_{\alpha, av}t}.$$
In summary, if $\alpha R_{\alpha, i}<0$ for all $i\in V$, there exists constants $c_1$ and $c_2$ such that
$$c_1e^{\alpha R_{\alpha, av}t}\leq R_{\alpha,i}-R_{\alpha, av}\leq c_2e^{\alpha R_{\alpha, av}t},$$
which implies the exponential convergence of the normalized $\alpha$-Ricci flow (\ref{normalized alpha flow}).
\qed\\

Using the normalized $\alpha$-Ricci flow (\ref{normalized alpha flow}),
the combinatorial $\alpha$-Calabi flow (\ref{definition of alpha Calabi flow})
and the modified $\alpha$-Calabi flow (\ref{modified alpha Calabi flow}),
we can give a characterization of the existence of constant $\alpha$-curvature metric.

\begin{theorem}\label{main alpha convergence theorem}
Suppose $(M, \mathcal{T}, \Phi)$ is a weighted triangulated surface with $\alpha\chi(M)\leq 0$. Then the existence of constant $\alpha$-curvature circle packing metric, the convergence of $\alpha$-Ricci flow (\ref{normalized alpha flow}), the convergence of $\alpha$-Calabi flow (\ref{definition of alpha Calabi flow}) and the convergence of modified $\alpha$-Calabi flow (\ref{modified alpha Calabi flow}) are mutually equivalent. Furthermore, if the solutions of the flows converge, then they all converge exponentially fast to a constant $R_{\alpha}$-curvature metric.
\end{theorem}
\textbf{Proof.} Assume there exists a constant curvature metric $r^*$, $u^*$ is the corresponding $u$-coordinate. Consider the $\alpha$-Ricci potential
\begin{equation}
F_\alpha(u)=\int_{u^*}^u\sum_{i=1}^N(K_i-R_{\alpha,av}r_i^\alpha)du_i.
\end{equation}
It has similar properties as stated in Lemma \ref{property of Ricci potential under general condition}.
Restricted to the hypersurface $\mathcal{U}=\{u\in\mathds{R}^N|\sum_i e^{\alpha u_i}=\sum_i e^{\alpha u_i(0)}\}$,
we have
$$\lim_{u\in\mathcal{U}, u\rightarrow \infty}F_\alpha(u)=+\infty$$
and $F_\alpha|_{\mathcal{U}}$ is also proper.
The rest proof is similar to that of Theorem \ref{global convergence of NCRF under general condition} and Theorem \ref{global convergence of CCF},
so we omit it here.\qed\\

We can also use the combinatorial $\alpha$-Ricci flow, combinatorial $\alpha$-Calabi flow and modified combinatorial $\alpha$-Calabi flow
to study the prescribing $\alpha$-curvature problem. Specifically, we have the following result.

\begin{theorem}\label{existence for prescribing alpha curvature problem}
Suppose $(M, \mathcal{T}, \Phi)$ is a weighted triangulated surface, $\alpha\in \mathds{R}$ is a
given real number, $\overline{R}_\alpha\in C(V)$
is a function defined on $M$ with $\alpha\overline{R}_{\alpha,i}\leq 0$ for $i=1, \cdots, N$.
Then $\overline{R}_\alpha$ is an admissible $\alpha$-curvature if and only if the
combinatorial $\alpha$-Ricci flow with target
\begin{equation*}
\begin{aligned}
\frac{dr_i}{dt}=(\overline{R}_{\alpha,i}-R_{\alpha,i})r_i,
\end{aligned}
\end{equation*}
exists for all time and converges,
if and only if the combinatorial $\alpha$-Calabi flow with target
\begin{equation*}
\frac{d r_i}{dt}=\Delta_\alpha (R_\alpha-\overline{R}_\alpha)_ir_i
\end{equation*}
exists for all time and converges, if and only if the modified combinatorial $\alpha$-Calabi flow with target
\begin{equation*}
\dot{u}=-\frac{1}{2}\nabla_u\overline{C}
\end{equation*}
exists for all time and converges, where $\overline{C}=\sum_i(K_i-\overline{R}_{\alpha,i}r_i^2)^2$.
\end{theorem}

\begin{remark}
For $\alpha=0$, the condition $\alpha\chi(M)\leq 0$ is always satisfied, so there is no restriction on $\chi(M)$ and the equivalences in Theorem \ref{main alpha convergence theorem} and Theorem \ref{existence for prescribing alpha curvature problem} are always valid.
For $\alpha=0$, the equivalent condition given by combinatorial Ricci flow is obtained in \cite{CL1}, and the equivalent condition given by combinatorial Calabi flow is obtained in \cite{Ge}.
\end{remark}

Now we consider the uniqueness and existence of constant $\alpha$-curvature metrics.
As the results are parallel to that of $R$-curvature, we just give the statements of the results and omit the proofs.

\begin{theorem}\label{uniqueness for alpha curvature}
Suppose $(M, \mathcal{T}, \Phi)$ is a weighted triangulated surface with $\alpha\chi(M)\leq 0$, then the constant $\alpha$-curvature metric is unique if it exists. Specificly, if $\alpha\chi(M)=0$, then there exists at most one constant $\alpha$-curvature metric up to scaling. If $\alpha\chi(M)<0$, then for any $c^*$, there exists at most one metric with $\alpha$-curvature $R_{\alpha,i}\equiv c^*$.
\end{theorem}

\begin{remark}
For $\alpha=0$, this is the result obtained by Thurston \cite{T1} and  Chow and Luo \cite{CL1}.
\end{remark}

\begin{theorem}\label{Thm-combtopo-condition-CCCPmetric}
Suppose $(M, \mathcal{T}, \Phi)$ is a weighted triangulated surface with $\alpha\chi(M)\leq 0$.
Then there exists a constant $\alpha$-curvature metric if and only if
there exists a circle packing metric $r^*$ such that, for any nonempty proper subset $I$ of $V$,
\begin{equation}\label{combtopo-condition of GX}
2\pi\chi(M)\frac{\sum_{i\in I}r_i^{*\alpha}}{\|r^*\|^{\alpha}_{\alpha}}>-\sum_{(e,v)\in Lk(I)}(\pi-\Phi(e))+2\pi\chi(F_I).
\end{equation}
\end{theorem}

\begin{theorem}\label{main existence thm for alpha-curvature}
Suppose $(M, \mathcal{T}, \Phi)$ is a weighted triangulated surface. If $\alpha\chi(M)<0$, then there exists a constant $\alpha$-curvature metric if and only if $\mathscr{Y}\cap\mathds{R}^N_{<0}\neq\phi$ when $\alpha>0$ and $\mathscr{Y}\cap\mathds{R}^N_{>0}\neq\phi$ when $\alpha<0$.
\end{theorem}

\begin{corollary}\label{Thm-uniqueness-CCCPmetric}
Suppose $(M, \mathcal{T}, \Phi)$ is a weighted triangulated surface with $\alpha\chi(M)>0$.
If there exists a constant $\alpha$-curvature metric,
then there exists a circle packing metric $r^*$ such that (\ref{combtopo-condition of GX}) is valid,
which implies that $\mathscr{Y}\cap\mathds{R}^N_{>0}\neq\phi$ in the case of $\alpha>0$ and $\chi(M)>0$,
and $\mathscr{Y}\cap\mathds{R}^N_{<0}\neq\phi$ in the case of $\alpha<0$ and $\chi(M)<0$.
\end{corollary}

\begin{theorem}\label{general existence theorem}
Suppose $(M, \mathcal{T}, \Phi)$ is a weighted triangulated surface. If there is a circle packing metric $r$ satisfying $R_{\alpha,i}\geq0$ for all $i\in V$, and
\begin{equation*}
-\sum_{(e,v)\in Lk(I)}(\pi-\Phi(e))+2\pi\chi(F_I)<0
\end{equation*}
for any nonempty subset $I$ of $V$, then there exists a non-negative constant $\alpha$-curvature metric.
\end{theorem}

\begin{remark}
Theorem \ref{existence theorem for negative curvature} and Theorem \ref{positive Euler number Thm}
are now special cases of Theorem \ref{alpha-converge thm for negative initial curvature} and Theorem \ref{general existence theorem}.
\end{remark}

\section{3-dimensional combinatorial Yamabe problem}\label{3-dimensional combinatorial Yamabe problem}

\subsection{The definition of combinatorial scalar curvature}

Suppose $M$ is a 3-dimensional compact manifold with a triangulation $\mathcal{T}=\{V,E,F,T\}$,
where the symbols $V,E,F,T$
represent the sets of vertices, edges, faces and tetrahedrons respectively.
A sphere packing metric is a map $r:V\rightarrow (0,+\infty)$ such that the length between
vertices $i$ and $j$ is $l_{ij}=r_{i}+r_{j}$ for each edge $\{i,j\}\in E$,
and the lengths $l_{ij},l_{ik},l_{il},l_{jk},l_{jl},l_{kl}$ determines a Euclidean
tetrahedron for each tetrahedron $\{i,j,k,l\}\in T$.
We can take sphere packing metrics as points in $\mathds{R}^N_{>0}$,
$N$ times Cartesian product of $(0,\infty)$, where $N=V^{\#}$ denotes the number of vertices.
It is pointed out \cite{G1} that a tetrahedron $\{i,j,k,l\}\in T$ generated by four positive radii $r_{i},r_{j},r_{k},r_{l}$
can be realized as a Euclidean tetrahedron if and only if
\begin{equation}\label{nondegeneracy condition}
Q_{ijkl}=\left(\frac{1}{r_{i}}+\frac{1}{r_{j}}+\frac{1}{r_{k}}+\frac{1}{r_{l}}\right)^2-
2\left(\frac{1}{r_{i}^2}+\frac{1}{r_{j}^2}+\frac{1}{r_{k}^2}+\frac{1}{r_{l}^2}\right)>0.
\end{equation}
Thus the space of admissible Euclidean sphere packing metrics is
$$\mathfrak{M}_{\mathcal{T}}=\left\{\;r\in\mathds{R}^N_{>0}\;\big|\;Q_{ijkl}>0, \;\forall \{i,j,k,l\}\in T\;\right\}.$$
Cooper and Rivin \cite{CR} called the tetrahedrons generated in this way conformal and
proved that a tetrahedron is conformal if and only if there exists a unique sphere
tangent to all of the edges of the tetrahedron. Moreover, the point of tangency with the
edge $\{i,j\}$ is of distance $r_i$ to $v_i$.
They further proved that $\mathfrak{M}_{\mathcal{T}}$ is a simply connected open subset
of $\mathds{R}^N_{>0}$, but not convex.

For a triangulated 3-manifold $(M, \mathcal{T})$ with sphere packing metric $r$, there is
also the notion of combinatorial scalar curvature.
Cooper and Rivin \cite{CR} defined the combinatorial scalar curvature $K_{i}$
at a vertex $i$ as angle deficit of solid angles
\begin{equation}\label{CR curvature}
K_{i}= 4\pi-\sum_{\{i,j,k,l\}\in T}\alpha_{ijkl},
\end{equation}
where $\alpha_{ijkl}$ is the solid angle at the vertex $i$ of the Euclidean tetrahedron $\{i,j,k,l\}\in T$
and the sum is taken over all tetrahedrons with $i$ as one of its vertices.
$K_i$ locally measures the difference between the volume growth rate of a small ball
centered at vertex $v_i$ in $M$ and a Euclidean ball of the same radius.
Cooper and Rivin's definition of combinatorial scalar curvature is motivated by the
fact that,
in the smooth case, the scalar curvature at a point $p$ locally measures the difference of
the volume growth rate of the geodesic ball with center $p$ to the Euclidean ball \cite{LP, Be}.

Similar to the two dimensional case, Cooper and Rivin's definition of combinatorial scalar
curvature $K_i$ is scaling invariant, which is not so satisfactory. The authors \cite{GX2}
once defined a new combinatorial scalar curvature as  $\frac{K_i}{r_i}$ on 3-dimensional
triangulated manifold $(M, \mathcal{T})$ with sphere packing metric $r$.
Motivated by the analysis in Section \ref{preliminaries of 2 dim}, we find that
it is more natural to define the combinatorial scalar curvature in the following way.

\begin{definition}\label{definition of comb scalar curv}
For a triangulated 3-manifold $(M, \mathcal{T})$ with a sphere packing metric $r$,
the combinatorial scalar curvature at the vertex $i$ is defined as
\begin{equation}\label{comb scalar curv}
R_i=\frac{K_i}{r_i^2},
\end{equation}
where $K_i$ is given by (\ref{CR curvature}).
\end{definition}

As $R_i$ differs from $K_i$ only by a factor $\frac{1}{r_i^2}$, $R_i$ still
locally measures the difference between the volume growth rate of a small ball
centered at vertex $v_i$ in $M$ and a Euclidean ball of the same radius.
Furthermore, according to the analysis in Section \ref{preliminaries of 2 dim},
$r_i^2$ is the analogue of the smooth Riemannian metric. If $\widetilde{r}_i^2=c r_i^2$
for some positive constant $c$, we have $\widetilde{R}_i=c^{-1}R_i$. This is
similar to the transformation of scalar curvature in the smooth case under scaling.

Analogous to the two dimensional case, we can define a measure on the vertices. As $r_i^2$
is the analogue of the Riemannian metric, we can take the measure as $\mu_i=r_i^3$, which
corresponds to the volume element. Then the total combinatorial scalar curvature is
\begin{equation}\label{total comb scalr curv}
\mathcal{S}=\int_M R d\mu=\sum_{i=1}^N R_i r_i^3=\sum_{i=1}^N K_i r_i.
\end{equation}
Note that $\mathcal{S}$ is just the functional introduced by Cooper and Rivin in \cite{CR}.
For the total combinatorial scalar curvature $\mathcal{S}$, we have the following
important property.

\begin{lemma}\label{property of Lambda}
(\cite{CR, Ri, G3})
Suppose $(M, \mathcal{T})$ is a triangulated 3-manifold with sphere packing metric $r$,
$\mathcal{S}$ is the total combinatorial scalar curvature. Then we have
\begin{equation}
\nabla_r\mathcal{S}=K.
\end{equation}
If we set
\begin{displaymath}
\Lambda=Hess_r\mathcal{S}=
\frac{\partial(K_{1},\cdots,K_{N})}{\partial(r_{1},\cdots,r_{N})}=
\left(
\begin{array}{ccccc}
 {\frac{\partial K_1}{\partial r_1}}& \cdot & \cdot & \cdot &  {\frac{\partial K_1}{\partial r_N}} \\
 \cdot & \cdot & \cdot & \cdot & \cdot \\
 \cdot & \cdot & \cdot & \cdot & \cdot \\
 \cdot & \cdot & \cdot & \cdot & \cdot \\
 {\frac{\partial K_N}{\partial r_1}}& \cdot & \cdot & \cdot &  {\frac{\partial K_N}{\partial r_N}}
\end{array}
\right),
\end{displaymath}
then $\Lambda$ is positive semi-definite with rank $N-1$ and
the kernel of $\Lambda$ is the linear space spanned by the vector $r$.
\end{lemma}

It should be emphasized that , as pointed out by Glickenstein \cite{G1}, the element $\frac{\partial K_i}{\partial r_j}$
for $i\sim j$ maybe negative, which is different from that of the two dimensional case.

By the definition of measure $\mu$, the total measure of $M$ is $\mu(M)=\sum_{i=1}^Nr_i^3$.
We will denote the total measure of $M$ by $V$ for simplicity in the following, if there is
no confusion.
The average combinatorial scalar curvature is
\begin{equation}\label{average comb scalar curv}
R_{av}=\frac{\int_M Rd\mu}{\int_Md\mu}=\frac{\mathcal{S}}{V}=\frac{\sum_{i=1}^NK_ir_i}{\sum_{i=1}^Nr_i^3}.
\end{equation}

\subsection{Combinatorial Yamabe problem in 3 dimension}

For the curvature $R_i$, it is natural to consider the corresponding constant curvature problem.
Suppose $R_i=\lambda, \forall i\in V,$ for some constant $\lambda$, then we have $K_i=\lambda r_i^2$,
which implies that $\lambda=\frac{\mathcal{S}}{V}$.
We can take $\frac{\mathcal{S}}{V}$ as a functional of the sphere packing metric $r$. However, this functional
is not scaling invariant in $r$. So we can modified the functional as $\frac{\mathcal{S}}{V^{1/3}}$. This
recalls us of the smooth Yamabe problem.

The classical Yamabe problem aims at solving the problem of existence of
constant scalar curvature metric on a closed manifold.
In order to study the constant scalar curvature problem on a closed Riemannian manifold $(M,g)$, Yamabe \cite{Ya} introduced
the so-called Yamabe functional $Q(g)$ and the Yamabe invariant $Y_{M,[g_0]}$, which are defined as
$$Q(g)=\frac{\int_M Rdv_g}{(\int_M dv_g)^{1-\frac{2}{n}}},$$
$$Y_{M,[g_0]}=\inf_{g\in[g_0]} Q(g),$$
where $[g_0]$ is the conformal class of Riemannian metric $g_0$. Trudinger \cite{Tr} and Aubin \cite{Au} made lots of
contributions to this problem, and Schoen \cite{Sc} finally gave the solution of the Yamabe Problem.
We refer the readers to \cite{LP} for this problem.

For piecewise flat manifolds, one can introduce similar functionals and invariants.
Champion, Glickenstein and Young \cite{CGY} studied Einstein-Hilbert-Regge functionals and
related invariants on triangulated manifolds with piecewise linear metrics. The authors \cite{GX2}
also introduced a type of combinatorial Yamabe functional and studied its properties.
To study the curvature defined by (\ref{comb scalar curv}), we can introduce the following
definition of combinatorial Yamabe functional and Yamabe invariant on triangulated 3-manifolds with
sphere packing metrics.

\begin{definition}\label{definition of comb Yamabe func}
Suppose $(M, \mathcal{T})$ is a triangulated 3-manifold with a fixed triangulation $\mathcal{T}$.
The combinatorial Yamabe functional is defined as
\begin{equation}
Q(r)=\frac{\mathcal{S}}{V^{\frac{1}{3}}}=\frac{\sum_{i=1}^NK_ir_i}{(\sum_{i=1}^Nr_i^3)^{1/3}}, \ \ r\in\mathfrak{M}_{\mathcal{T}}.
\end{equation}
The combinatorial Yamabe invariant with respect to $\mathcal{T}$ is defined as
$$Y_{M,\mathcal{T}}=\inf_{r\in\mathfrak{M}_{\mathcal{T}}} Q(r).$$
\end{definition}

The admissible sphere packing metric space $\mathfrak{M}_{\mathcal{T}}$ for a given triangulated manifold $(M,\mathcal{T})$
is an analogue of the conformal class $[g_0]$ of a Riemannian manifold $(M, g_0)$, as every admissible sphere packing
metric could be taken to be conformal to the metric with all $r_i=1$. We call $\mathfrak{M}_{\mathcal{T}}$
the combinatorial conformal class for $(M,\mathcal{T})$. It is uniquely determined by the triangulation $\mathcal{T}$
of $M$. $Y_{M,\mathcal{T}}$ is referred to as the Yamabe constant for $(M,\mathcal{T})$.

Note that the combinatorial Yamabe invariant $Y_{M,\mathcal{T}}$ is an invariant of the conformal class $\mathfrak{M}_{\mathcal{T}}$
and is well defined, as we have
\begin{equation}\label{Q<|K|}
|Q(r)|=\left|\frac{\mathcal{S}}{V^{1/3}}\right|=\left|\frac{\Sigma K_ir_i}{\left(\Sigma r_i^3\right)^{\frac{1}{3}}}\right|\leq \left(\Sigma K_i^{3/2}\right)^{\frac{2}{3}}=\|K\|_{3/2},
\end{equation}
and $K_i$ satisfies $(4-2d)\pi\leq K_i <4\pi$,
where $d$ is the maximal degree and $d\leq E^{\#}$.
If the equality in (\ref{Q<|K|}) is achieved, the corresponding sphere packing metric
must be a constant curvature metric.

By direct computation, we have
\begin{equation}\label{gradient of Q}
\nabla_{r_i}Q=\frac{1}{V^{1/3}}(K_i-R_{av}r_i^2),
\end{equation}
which implies that $r$ is a constant combinatorial scalar curvature metric if and only if it is
a critical point of combinatorial Yamabe functional $Q(r)$.

Analogue to the smooth Yamabe problem, we can raise the following combinatorial Yamabe problem
on 3-dimensional triangulated manifold.\\
~\\
\textbf{The Combinatorial Yamabe Problem.}
Given a 3-dimensional manifold $M$ with a triangulation $\mathcal{T}$,
find a sphere packing metric with constant combinatorial scalar curvature in the combinatorial
conformal class $\mathfrak{M}_{\mathcal{T}}$.\\

We could further consider finding a suitable triangulation for $M$ which admits a constant combinatorial scalar curvature metric.

It is easy to see that the Platonic solids with tetrahedral cells all admit constant combinatorial scalar curvature metric,
including the 5-cell, the 16-cell, the 600-cell, etc. In these cases, the constant combinatorial scalar curvature metrics arise from
symmetry and taking the radii equal.

\subsection{Combinatorial Yamabe flow}

In 1980s, Hamilton \cite{Ham1, Ham3} proposed the Yamabe flow to the Yamabe problem. For a closed n-dimensional
Riemannian manifold $(M^n, g)$ with $n\geq 3$, the Yamabe flow is defined to be
\begin{equation}\label{Yamabe flow}
\frac{\partial}{\partial t}g_{ij}=-Rg_{ij}
\end{equation}
with normalization
\begin{equation}\label{normalized Yamabe flow}
\frac{\partial}{\partial t}g_{ij}=(r-R)g_{ij},
\end{equation}
where $R$ is the scalar curvature of $g$ and $r$ is the average of the scalar curvature.
Along the normalized Yamabe flow, the volume is invariant, the total scalar curvature is decreased, and
the scalar curvature evolves according to
\begin{equation}\label{evolution of R along Yamabe flow}
\frac{\partial R}{\partial t}=(n-1)\Delta R+R(R-r).
\end{equation}
It is proved by Hamilton \cite{Ham3} that the solution to the Yamabe flow (\ref{normalized Yamabe flow})
exists for all time and the solution converges exponentially fast to a metric with constant scalar
curvature if $R<0$ initially. Ye \cite{Ye} then proved that, if the initial metric is locally conformally
flat, then the solution to the normalized Yamabe flow (\ref{normalized Yamabe flow})
converges to a metric of constant scalar curvature.
Brendle \cite{Br1} proved that the solution of (\ref{normalized Yamabe flow}) also converges to
a metric of constant curvature if $3\leq n\leq 5$. He \cite{Br2} further handled the case of $n\geq 6$
and got some convergence results.

In the combinatorial case, Luo \cite{L1} first introduced the combinatorial Yamabe flow on
surfaces and Glickenstein \cite{G1, G2} introduced the combinatorial Yamabe flow on 3-dimensional
manifolds and studied its related properties.
Luo \cite{L2} also introduced a combinatorial curvature flow for piecewise constant
curvature metrics on compact triangulated 3-manifolds with boundary
consisting of surfaces of negative Euler characteristic.
To study the constant curvature problem of $R_i$, we
introduce the following combinatorial Yamabe flow.

\begin{definition}
Given a 3-dimensional triangulated manifold $(M, \mathcal{T})$ with sphere packing metric $r$,
the combinatorial Yamabe flow is defined to be
\begin{equation}\label{comb Yamabe flow}
\frac{dg_i}{dt}=-R_ig_i,
\end{equation}
with normalization
\begin{equation}\label{normalized comb Yamabe flow}
\frac{dg_i}{dt}=(R_{av}-R_i)g_i,
\end{equation}
where $g_i=r_i^2$ and $R_{av}$ is the average of the combinatorial scalar curvature given by (\ref{average comb scalar curv}).
\end{definition}

Following the 2-dimensional case, it is easy to check that the solutions of (\ref{comb Yamabe flow})
and (\ref{normalized comb Yamabe flow}) can be transformed to each other by a scaling procedure.
And it is easy to see that, if the solution of (\ref{normalized comb Yamabe flow}) converges to a sphere packing
metric $r(+\infty)$, then $r(+\infty)$ is a metric with constant combinatorial scalar curvature.
Analogous to the smooth Yamabe flow, we have the following properties.

\begin{proposition}
Suppose $r(t)$ is a solution of (\ref{normalized comb Yamabe flow}) on a triangulated 3-manifold
$(M, \mathcal{T})$. Along the flow (\ref{normalized comb Yamabe flow}), the total measure
$\mu(M)=V=\sum_{i=1}^Nr_i^3$ is invariant and the total combinatorial scalar curvature $\mathcal{S}$
is decreased.
\end{proposition}

\textbf{Proof.}
The normalized combinatorial Yamabe flow could be written as
$$\frac{dr_i}{dt}=\frac{1}{2}(R_{av}-R_i)r_i.$$
Using this equation, we have
$$\frac{dV}{dt}=3\sum_{i=1}^Nr_i^2\frac{dr_i}{dt}=\frac{3}{2}\sum_{i=1}^Nr_i^3(R_{av}-R_i)=0.$$
For the total combinatorial scalar curvature $\mathcal{S}$, we have
\begin{equation}\label{derivative of total comb scalar curv}
\begin{aligned}
\frac{d\mathcal{S}}{dt}
=&\sum_{i=1}^N\frac{dK_i}{dt}r_i+\sum_{i=1}^NK_i\frac{dr_i}{dt}\\
=&\frac{1}{2}\sum_{i,j=1}^Nr_i\frac{\partial K_i}{\partial r_j}(R_{av}-R_j)r_j+\frac{1}{2}\sum_{i=1}^NK_ir_i(R_{av}-R_i)\\
=&-\frac{1}{2}\sum_{i=1}^N\frac{K_i(K_i-R_{av}r_i^2)}{r_i}\\
=&-\frac{1}{2}\sum_{i=1}^N\frac{(K_i-R_{av}r_i^2)^2}{r_i},
\end{aligned}
\end{equation}
which implies that $\mathcal{S}$ is decreased along the flow (\ref{normalized comb Yamabe flow}).
Note that Lemma \ref{property of Lambda} is used in the third step. \qed\\

\begin{remark}
As the total measure of $M$ is invariant along (\ref{normalized comb Yamabe flow}), we will focus on the
properties of  (\ref{normalized comb Yamabe flow}) in the following. Furthermore, we will assume
$r(0)\in \mathbb{S}^{N-1}\triangleq\{r\in \mathds{R}^N; ||r||_3=(\sum r_i^3)^{\frac{1}{3}}=1\}$ in the following.
\end{remark}

By direct calculations, we find that the curvature $R_i$ evolves according to the following equation along the
normalized combinatorial Yamabe flow
(\ref{normalized comb Yamabe flow})
\begin{equation}\label{evolution of R_i in 3 dim a}
\frac{dR_i}{dt}=-\frac{1}{2r_i^2}\sum_{j=1}^N \frac{\partial K_i}{\partial r_j}r_jR_j+R_i(R_{i}-R_{av}).
\end{equation}
If we define the Laplacian as
\begin{equation}\label{definition of Laplacian in 3 dim}
\Delta f_i=-\frac{1}{r_i^2}\sum_{j=1}^N \frac{\partial K_i}{\partial r_j}r_jf_j
=\frac{1}{r_i^2}\sum_{j\sim i} (-\frac{\partial K_i}{\partial r_j}r_j)(f_j-f_i)
\end{equation}
for $f\in C(V)$, then the equation (\ref{evolution of R_i in 3 dim a}) could be written as
\begin{equation}\label{evolution of R_i along CYF}
\frac{dR_i}{dt}=\frac{1}{2}\Delta R_i+R_i(R_{i}-R_{av}),
\end{equation}
which has almost the same form as the evolution equation (\ref{evolution of R along Yamabe flow}) of scalar curvature
along the Yamabe flow (\ref{normalized Yamabe flow}) on three dimensional smooth manifolds.
The Laplacian defined by (\ref{definition of Laplacian in 3 dim}) satisfies $\int_M \Delta f d\mu=0$ and $\Delta c=0$ for any
$f\in C(V)$ and constant $c\in \mathds{R}$. Furthermore, by the calculations in \cite{G1}, we have
\begin{equation}\label{Laplacian using dual in 3 dim}
\Delta f_i=\frac{1}{r_i^3}\sum_{j\sim i}\frac{l_{ij}^*}{l_{ij}}(f_j-f_i),
\end{equation}
where $l_{ij}^*$ is the area dual to the edge $\{i, j\}$.
Note that this form of Laplacian is very similar to
Hirani's definition \cite{Hi} of Laplace-Beltrami operator
$$\Delta f_i=\frac{1}{V_i^*}\sum_{j\sim i}\frac{l_{ij}^*}{l_{ij}}(f_j-f_i),$$
except the first factor. It should be mentioned that $r_i^3$ is a type of volume.

Though the definition (\ref{definition of Laplacian in 3 dim}) of Laplacian has lots of good properties,
it is not a Laplacian on graphs in the usual sense, as $l_{ij}^*$ may be negative,
which makes the maximum principle for (\ref{evolution of R_i along CYF}) not so good.
This was founded by Glickenstein and he studied the properties of such Laplacian in \cite{G1,G2}.

To study the long time behavior of (\ref{normalized comb Yamabe flow}), we need to classify the solutions of the flow.

\begin{definition}
A solution $r(t)$ of the combinatorial Yamabe flow (\ref{normalized comb Yamabe flow})
is nonsingular if the solution $r(t)$ exists for
$t\in [0, +\infty)$ and $\{r(t)\}\subset\subset \mathfrak{M}_{\mathcal{T}}\cap \mathbb{S}^{N-1}$.
\end{definition}

In fact, the condition $\{r(t)\}\subset\subset \mathfrak{M}_{\mathcal{T}}\cap \mathbb{S}^{N-1}$
ensures the long time existence of the flow (\ref{normalized comb Yamabe flow}).
Furthermore, we have the following property for nonsingular solutions.

\begin{theorem}\label{existence of const curv metr in 3d}
If there exists a nonsingular solution for the flow (\ref{normalized comb Yamabe flow}), then there
exists at least one sphere packing metric with constant combinatorial scalar curvature $R_i$ on $(M, \mathcal{T})$.
\end{theorem}

\textbf{Proof.}
By (\ref{derivative of total comb scalar curv}), we know that $Q(r)$ is decreasing along the flow
(\ref{normalized comb Yamabe flow}). As $Q(r)$ is uniformly bounded by (\ref{Q<|K|}), the limit
$\lim_{t\rightarrow +\infty}Q(r(t))$ exists. Then there exists a sequence $t_n\uparrow +\infty$ such
that
\begin{equation}\label{derivative of Q}
(Q(r))'(t_n)=-\frac{1}{2}\sum_{i=1}^N\frac{(K_i-R_{av}r_i^2)^2}{r_i}\rightarrow 0.
\end{equation}
As $\{r(t)\}\subset\subset \mathfrak{M}_{\mathcal{T}}\cap \mathbb{S}^{N-1}$, there exists
a subsequence, denoted as $r_n$, of $r(t_n)$ such that
$r_n\rightarrow r^*\in \mathfrak{M}_{\mathcal{T}}\cap \mathbb{S}^{N-1}$.
Then (\ref{derivative of Q}) implies that $r^*$ satisfies $K_i^*=R_{av}(r_i^*)^2$ and
$r^*$ is a sphere packing metric with constant combinatorial scalar curvature $R_i$.\qed\\

In fact, under some conditions, we find that the sphere packing metrics with constant combinatorial
scalar curvature are isolated on $\mathbb{S}^{N-1}$.

\begin{theorem}\label{isolation of cons comb scal curv metric in 3d}
The sphere packing metrics with nonpositive constant combinatorial scalar curvature
are isolated in $\mathfrak{M}_{\mathcal{T}}\cap \mathbb{S}^{N-1}$.
\end{theorem}
\textbf{Proof.}
Define
\begin{equation*}
\begin{aligned}
G:\mathfrak{M}_{\mathcal{T}}&\rightarrow \mathds{R}^{N}\\
r&\mapsto \left(\frac{1}{r_1}(K_1-\lambda r_1^2), \cdots, \frac{1}{r_N}(K_N-\lambda r_N^2)\right),
\end{aligned}
\end{equation*}
where $\lambda=\frac{\mathcal{S}}{V}$.
It is easy to see that the zero point of $G$ corresponds to the metric with
constant combinatorial scalar curvature. By direct calculations, the Jacobian of $G$ at the constant
curvature metric is
\begin{equation}
DG=\left(\Lambda-2\lambda\Sigma^{\frac{1}{4}}\left(I-\frac{r^{3/2}\cdot (r^{3/2})^T }{V}\right)\Sigma^{\frac{1}{4}}\right)\Sigma^{-\frac{1}{2}},
\end{equation}
where $r^{3/2}=(r_1^{3/2},\cdots,r_N^{3/2})^T$.
If $r^*$ is a sphere packing metric with nonpositive constant combinatorial scalar
curvature, then $\lambda\leq0$ and hence $DG$ is a positive semi-definite matrix with rank $N-1$ and kernel $\{cr^2|c\in \mathds{R}\}$,
which is normal to $\mathbb{S}^{N-1}$. Restricted to $\mathbb{S}^{N-1}$,
$DG$ is positive definite and then nondegenerate, which implies that the
zero points of $G$ with nonpositive curvature are isolated in $\mathfrak{M}_{\mathcal{T}}\cap \mathbb{S}^{N-1}$. \qed\\

Suppose the solution of the flow (\ref{normalized comb Yamabe flow})
exists for $t\in [0,T)$ with $0<T\leq +\infty$, then we have
$\{r(t)\}\subseteq \mathfrak{M}_{\mathcal{T}}\cap \mathbb{S}^{N-1}$. However, if
$$\overline{\{r(t)\}}\cap \partial(\mathfrak{M}_{\mathcal{T}}\cap \mathbb{S}^{N-1})\neq\emptyset,$$
the flow will raise singularities. They could be separated into two types as follows.
\begin{description}
  \item[Essential Singularity] There is a vertex $i\in V$ such that there exists a sequence of time $t_n\uparrow T$
such that $r_i(t_n)\rightarrow 0$ as $n\rightarrow +\infty$;
  \item[Removable Singularity] There exists a sequence of time $t_n\uparrow T$ such that
$r_i(t_n)\subset\subset \mathds{R}_{>0}$ for all vertices $i\in V$, but there exists a tetrahedron
$\{ijkl\}\in \mathcal{T}$ such that $Q_{ijkl}(t_n)\rightarrow 0$ as $n\rightarrow +\infty$.
\end{description}

It seems that the following conjectures are likely to hold for
the normalized combinatorial Yamabe flow (\ref{normalized comb Yamabe flow})  on 3-dimensional manifolds.
\conjecture The normalized combinatorial Yamabe flow (\ref{normalized comb Yamabe flow}) will not develop essential singularity in finite time.
\conjecture \label{conjecture2}If no singularity develops along the normalized flow (\ref{normalized comb Yamabe flow}),
the solution converges to a sphere packing metric with constant combinatorial scalar curvature as time approaches infinity.\\

It is interesting to note that Glickenstein \cite{G2} made a small amount of progress on Conjecture 1.
He proved that, for some special class of complexes, the flow
he introduced would not develop  singularities in finite time.

We have mentioned above that the existence of sphere packing metric with constant combinatorial
scalar curvature is necessary for the convergence of the normalized combinatorial Yamabe flow (\ref{normalized comb Yamabe flow}).
In fact, it is almost sufficient and we have the following result.

\begin{theorem}\label{convergence of CYF under existence}
Suppose $r^*$ is a sphere packing metric on $(M, \mathcal{T})$ with nonpositive constant combinatorial
scalar curvature. If $||r(0)-r^*||^2$ is small enough, the solution of the normalized combinatorial
Yamabe flow (\ref{normalized comb Yamabe flow}) exists for all time and converges to $r^*$.
\end{theorem}

\textbf{Proof.}
We can rewrite the normalized combinatorial Yamabe flow (\ref{normalized comb Yamabe flow}) as
$$\frac{dr_i}{dt}=\frac{1}{2}(R_{av}-R_i)r_i.$$
Set $\Gamma_i(r)=\frac{1}{2}(R_{av}-R_i)r_i$, then by the calculations in the proof of Theorem
\ref{isolation of cons comb scal curv metric in 3d}, we have
$$D\Gamma|_{r^*}=-\frac{1}{2}Z\cdot \left(
                           \begin{array}{ccc}
                             \frac{1}{r_1} &   &   \\
                               & \ddots &   \\
                               &   & \frac{1}{r_N} \\
                           \end{array}
                         \right),$$
where
$$Z=\Lambda-2R_{av}\left(
                           \begin{array}{ccc}
                             r_1^{1/2} &   &   \\
                               & \ddots &   \\
                               &   & r_N^{1/2} \\
                           \end{array}
                         \right)
\left(I-\frac{r^{3/2}\cdot (r^{3/2})^T }{V}\right)
\left(
                           \begin{array}{ccc}
                             r_1^{1/2} &   &   \\
                               & \ddots &   \\
                               &   & r_N^{1/2} \\
                           \end{array}
                         \right).$$
If $R_{av}\leq 0$, then $D\Gamma|_{r^*}$ is a matrix with rank $N-1$ and kernel $\{cr^2|c\in \mathds{R}\}$.
Furthermore, the nonzero eigenvalues of $D\Gamma|_{r^*}$ are all negative.
Note that, along the normalized flow (\ref{normalized comb Yamabe flow}), $||r||_3$
is invariant and thus the kernel $\{cr^2|c\in \mathds{R}\}$ is transversal to the flow. This implies that
$D\Gamma|_{r^*}$ is negative definite on $\mathbb{S}^{N-1}$ and $r^*$ is
a local attractor of the normalized combinatorial Yamabe flow (\ref{normalized comb Yamabe flow}).
Then the conclusion follows from the Lyapunov Stability Theorem. \qed\\

This theorem has the following interesting corollary.

\begin{corollary}\label{convergence of CYF under S(0)<0}
Given a 3-dimensional triangulated manifold $(M^3,\mathcal{T})$.
If the initial total combinatorial scalar curvature functional $\mathcal{S}(0)\leq 0$,
and no singularity develops along the normalized combinatorial Yamabe flow (\ref{normalized comb Yamabe flow}),
then the solution converges to a sphere packing metric $r^*$ with nonpositive constant combinatorial
scalar curvature.
\end{corollary}
\textbf{Proof.}
By the proof of Theorem \ref{existence of const curv metr in 3d},
if the solution $r(t)$ of (\ref{normalized comb Yamabe flow}) is nonsingular,
there exists a sequence $r_n\rightarrow r^*$, where $r^*$ is a sphere packing
metric with constant combinatorial scalar curvature. As $\mathcal{S}$ is decreasing
along the flow (\ref{normalized comb Yamabe flow}) and $\mathcal{S}(0)\leq 0$,
$r^*$ has nonpositive combinatorial scalar curvature. Then Theorem \ref{convergence of CYF under existence}
implies the conclusion of the corollary. \qed

\begin{remark}
Especially, if $R_i(0)\leq 0$, Corollary \ref{convergence of CYF under S(0)<0} is still valid.
\end{remark}

\section{Some questions}\label{unsolved problems}
There are several questions which we find interesting relating to the results in the paper.

\begin{enumerate}
  \item  From the results in Section \ref{2-d comb RF},
we see that the surfaces with positive Euler characteristic are
particular for the  constant curvature problem on surfaces.
Example \ref{example triangluate sphere} shows that combinatorial Ricci flow (\ref{equation of CCF introduction})
may not converge to a constant $R$-curvature metric,
while Example \ref{example-not unique cccpm} shows that constant curvature metric may not be unique.
If we want to approximate a smooth geometric object by corresponding discrete objects,
the quantities on these objects should obey similar laws or exhibit similar properties.
However, it has been proved \cite{Ham1, CH1} that the Ricci flow on surfaces converges to a constant curvature metric for
any initial metric.
The differences between discrete case and the smooth case are very interesting. It deserves deeper studying.
Maybe these differences are caused by the triangulation. The tetrahedron triangulation used in Example \ref{example triangluate sphere}
and \ref{example-not unique cccpm}
maybe too rough to approximate the sphere.
We expect these differences will disappear as the triangulation becomes finer and finer.
We want to know whether we can find a triangulation for the sphere such that the constant curvature metric uniquely exists.
We further want to know how to evolve discrete curvature $R_i=K_i/r_i^2$ along discrete curvature flows to constant curvature or, more generally,
evolve discrete $\alpha$-curvature to constant curvature when $\alpha\chi(M)>0$.

  \item Whether there are similar topological and combinatorial obstructions, similar to Thurston's criterion (\ref{Thuston's condition}),
for the existence of constant $\alpha$-curvature metric?
For the negative constant curvature metric, we derived an existence result (Theorem \ref{existence theorem for negative curvature})
by the discrete maximum principe, which is different from the other results.
We want to further develop these methods to derive more existence results.
We want to know whether we can get similar convergence results as that of the smooth surface Ricci flow and
then give a discrete uniformization theorem for $R_i$.
We believe that the results in \cite{BHLLMY} will play an important role in the procedure if it is feasible.
We also want to know the relationship between our conditions and Thurston's criterion (\ref{Thuston's condition}).

  \item The first positive eigenvalue of discrete Laplace operator
  plays an important role in the proof of the main results in two dimension.
 It is interesting to estimate the first positive eigenvalue of discrete Laplace operator
  and use it to derive some existence results.

  \item Investigate the prescribing curvature problem when the target curvature $\overline{R}_i>0$ or, more generally, $\alpha \overline{R}_{\alpha, i}>0$ for at least one vertex $i$.
  \item Study the singularities of three dimensional combinatorial Yamabe flow (\ref{normalized comb Yamabe flow introduction}), and find topological and combinatorial obstructions for the existence of sphere packing metric with constant combinatorial scalar curvature.
  \item Study the discrete Yamabe functional (\ref{definition of comb Yamabe func}) and solve the combinatorial Yamabe problem in three dimension.
  \item Study the combinatorial Yamebe problem in higher dimensions.
\end{enumerate}

\textbf{Acknowledgements}\\[8pt]
The first author would like to show his greatest respect to Professor Gang Tian who brought him to
the area of combinatorial curvature flows.
The research of the second author is partially supported by National Natural Science Foundation of China under
grant no. 11301402 and 11301399.
He would also like to thank Professor Guofang Wang for the invitation to the Institute
of Mathematics of the University of Freiburg and for his encouragement and many
useful conversations during the work.
Both authors would also like to thank Dr. Dongfang Li for many helpful conversations.

(Huabin Ge)
Department of Mathematics, Beijing Jiaotong University, Beijing 100044, P.R. China

E-mail: hbge@bjtu.edu.cn\\[2pt]

(Xu Xu) School of Mathematics and Statistics, Wuhan University, Wuhan 430072, P.R. China

E-mail: xuxu2@whu.edu.cn\\[2pt]

\end{document}